\documentclass[review, 3p]{elsarticle}

\usepackage{hyperref}
\usepackage{units}
\usepackage{amsmath}
\usepackage{amsmath,amsfonts,amsthm} 
\usepackage{units}
\usepackage{multirow}
\usepackage{graphicx}
\usepackage{float}
\usepackage{subcaption}
\usepackage{array}
\usepackage{enumitem}
\usepackage{epstopdf}
\usepackage{soul}

\newcommand{\mathtext}[1]{\text{#1}}

\newcommand{\defineSet}[2]{\{#1 \hspace{2pt}:\hspace{2pt} #2 \}}

\newcommand{\partialDt}[1]{\frac{\partial#1}{\partial t}}
\newcommand{\Divergence}[1]{\nabla \cdot \boldsymbol{#1}}
\newcommand{\GradientScalar}[1]{\nabla #1 }
\newcommand{\GradientVector}[1]{\nabla \boldsymbol{#1} }
\newcommand{\areaIntegral}[2]{\int\limits_{#1}^{}#2\mathtext{d}\Omega}
\newcommand{\lineIntegral}[2]{\oint\limits_{#1}^{}#2\mathtext{d}\varGamma}

\newcommand{\DirichletBC}{\varGamma_D}
\newcommand{\NeumannBC}{\varGamma_N}
\newcommand{\Real}{\mathbb{R}}
\newcommand{\xDiscreteScalar}[1]{#1_h}
\newcommand{\xtDiscreteScalar}[2]{#1_h^{#2}}
\newcommand{\ContVector}[1]{\boldsymbol{#1}}
\newcommand{\xDiscreteVector}[1]{\boldsymbol{#1}_h}

\newcommand{\xtDiscreteVector}[2]{\boldsymbol{#1}_h^{#2}}

\newcommand{\brac}[1]{\left( {#1} \right)}
\newcommand{\InDomain}[1]{\hspace{10pt} \text{in} \hspace{10pt} #1}
\newcommand{\OnBoundary}[1]{\hspace{10pt} \text{on} \hspace{10pt} #1}

\makeatletter
\def\old@comma{,}
 \catcode`\,=13
 \def,{%
   \ifmmode%
     \old@comma\discretionary{}{}{}%
   \else%
     \old@comma%
   \fi%
 }
 \makeatother

\journal{Elsevier}

\begin{document}
\begin{frontmatter}
\title{Constrained particle-mesh projections in a hybridized discontinuous Galerkin framework with applications to advection-dominated flows}
%
%
\author[a]{Jakob M. Maljaars \corref{cor1}}
\ead{j.m.maljaars@tudelft.nl}
\author[a]{Robert Jan Labeur}
\ead{r.j.labeur@tudelft.nl}
\author[b]{Nathaniel Trask \fnref{fn1}}
\ead{natrask@sandia.gov}
\author[c]{Deborah Sulsky}
\ead{sulsky@math.unm.edu}
\address[a]{Environmental~Fluid~Mechanics, Faculty~of~Civil~Engineering~and~Geosciences, Delft~University~of~Technology, Stevinweg~1, 2600~GA Delft, The~Netherlands}
\address[b]{Center~for~Computing~Research, Sandia~National~Laboratories, Albuquerque, NM~87185-1320, USA}
\address[c]{Department~of~Mathematics~and~Statistics, The~University~of~New~Mexico, Albuquerque, NM~87131, USA}
\cortext[cor1]{Corresponding author}
\fntext[fn1]{Sandia National Laboratories is a multimission laboratory managed and operated by National Technology and Engineering Solutions of Sandia, LLC, a wholly owned subsidiary of Honeywell International, Inc., for the U.S. Department of Energy's National Nuclear Security Administration under contract DE-NA0003525.}
\begin{abstract}
By combining concepts from particle-in-cell (PIC) and hybridized discontinuous Galerkin (HDG) methods, we present a particle-mesh scheme which allows for diffusion-free advection, satisfies mass and momentum conservation principles in a local sense, and allows the extension to high-order spatial accuracy.
To achieve this, we propose a novel particle-mesh projection operator required for the exchange of information between the particles and the mesh. Key is to cast these projections as a PDE-constrained $\ell^2$-optimization problem to allow the advective field naturally located on Lagrangian particles to be expressed as a mesh quantity. By expressing the control variable in terms of single-valued functions at cell interfaces, 
this optimization problem seamlessly fits in a HDG framework. Owing to this framework, the resulting scheme can be implemented efficiently via static condensation. The performance of the scheme is demonstrated by means of various numerical examples for the linear advection-diffusion equation and the incompressible Navier-Stokes equations. The results show that optimal spatial accuracy can be achieved, and given the particular time-stepping strategy, second-order time accuracy is confirmed. The robustness of the scheme is illustrated by considering benchmarks for advection of discontinuous fields and the Taylor-Green vortex instability in the high Reynolds number regime. 
\end{abstract}
\end{frontmatter}
\section{Introduction}
Flow and transport phenomena in which advection plays a dominant role, and that can be modeled as incompressible, arise in many engineering applications, such as the advection and mixing of pollutants, turbulent flows, and multiphase flows. 
Accurately simulating such problems requires a discretization of the underlying mass and momentum conservation laws having minimal artificial dissipation relative to the physical damping, as well as imposing the incompressibility constraint. 
Whereas Lagrangian particle-based methods offer a diffusion-free treatment of advective phenomena,  pure, Lagrangian particle motion consistent with the incompressibility constraint is difficult, if not impossible to achieve. On the other hand, Eulerian mesh-based methods are amenable to imposition of incompressibility and possess global or local (i.e. cellwise) conservation properties. However, in the advection dominated regime, the stabilization of mesh-based methods is non-trivial and typically introduces artificial diffusion. By combining concepts from Lagrangian particle-based and Eulerian mesh-based methods, this paper formulates a novel particle-mesh scheme that reconciles a diffusion-free discretization of the advective part of the problem with an efficient discretization of the constitutive part, while rigorously satisfying local conservation principles. Furthermore, the method extends to arbitrarily high-order spatial accuracy on structured or unstructured meshes.

Tracing back to the 1960s, the particle-in-cell (PIC) method \cite{Evans1957} was the first to consider Lagrangian particles in the advective part of a problem, while conveniently using a mesh to account for the dynamic interaction between the particles. The coupling between the particles and the mesh is established by means of two auxiliary projection steps: a particle-mesh projection for reconstructing field data from the scattered particles, and a mesh-particle projection to update the properties of the Lagrangian particles from the solution on the mesh.
Conservation in PIC and related methods, such as the material-point method (MPM) \cite{Sulsky1994}, is generally examined in the context of global conservation \cite{Brackbill1986,Burgess1992,Love2006,Jiang2017}. The salient principle is that the particle discretization of the solution and the mesh-based discretization of the solution should contain the same total mass, momentum and energy, as these are two representations of the same material. 
Indeed, exact conservation of these quantities is possible when using basic, low-order methods in which particles are used as moving point masses or volumes \cite{Love2006}. 
However, this does not yet insure that the particle fields are locally consistent with the mesh description of the flow. For example, the mass or momentum represented by the particles within a discretization mesh cell does not typically match the mass and momentum associated with the cell.
Extension to higher-order accuracy is possible by considering the particles as moving sampling points of the material at the mesh \cite{Wallstedt2011,Edwards2012,Sulsky2016,Maljaars2017}. While being consistent, this compromises however exact conservation \cite{Edwards2012,Maljaars2017}. 
This trade-off between conservation and accuracy is also observed in other particle methods such as SPH \cite{Monaghan2005,Dilts2003}.

In this work, as in related previous work \cite{Maljaars2017}, 
we frame our particle-mesh scheme as an operator splitting method and use particles as moving sampling points for reconstructing fields at the mesh. By extending the concepts explored in \cite{Maljaars2017}, this paper formulates particle-mesh projection operators which provide high-order accuracy, while simultaneously satisfying rigorous global and local conservation principles on the mesh. 
Not surprisingly, a minimal requirement for obtaining the latter is to employ a mesh-based framework exhibiting local conservation properties. This requirement is met by many meshed-based discretizations, and in particular is met by discontinuous Galerkin (DG) methods (see, e.g., \cite{Cockburn2000,Cockburn2005} and references) and hybridized discontinuous Galerkin (HDG) methods (see, e.g., \cite{Labeur2007,Nguyen2009,Egger2010,Labeur2012,Lehrenfeld2016,Rhebergen2017}).  
Typical to (H)DG is to pose the problem in terms of cellwise balances augmented with properly defined numerical fluxes at cell interfaces. These numerical fluxes are usually expressed in terms of the element field variables, and require upwinding (e.g., \cite{Labeur2007,Egger2010,Labeur2012,Lehrenfeld2016,Rhebergen2017}) or other carefully designed formulations (e.g., \cite{Nguyen2009,Nguyen2011}) in order to obtain stable solutions in advection dominated regimes. 
What distinguishes HDG from DG is that the former employs
single-valued interface functions which are defined on cell facets only. Typically, the interface flux in HDG is partly formulated in terms of these interface fields, thereby providing the conditions to enforce inter-element flux continuity in a weak sense.

Essential to our novel particle-mesh projection is the reconstruction of HDG interface fields from the advection of Lagrangian particles. In this way, we inherit the diffusion-free properties of a particle-based method while retaining the local conservation properties of a HDG method. The key idea in formulating the particle-mesh projection is that the particle motion has to satisfy a mesh-based advection operator. This advection operator is used as a constraint to augment the local $\ell^2$-projections used in \cite{Maljaars2017}, resulting in a PDE-constrained minimization problem for the particle-mesh projections. By expressing the control variable in terms of single-valued interface functions, this approach seamlessly fits in a HDG framework. The HDG framework provides an additional benefit; it enables an efficient implementation procedure via static condensation.

We will present our particle-mesh scheme for the advection-diffusion equation and the incompressible Navier-Stokes equations, where the inclusion of the diffusion terms and/or the incompressibility constraint is formulated as an operator splitting scheme. Considering the simpler case of advection-diffusion 
allows analysis of our scheme, before using our approach for more demanding problems, including the challenging case of vortex instability in high Reynolds number flows.
Other potential applications, beyond the scope of this paper, include the sharp representation, using particles, of fluid-fluid interfaces in multiphase flows.

The remainder of this work is structured as follows. In Section~\ref{sec:problem-definitions} the problems of interest are formulated and a spatial-temporal operator splitting technique is introduced. Formulating the PDE-constrained particle-mesh projection, and proving consistency and conservation constitute the main part of Section~\ref{sec: Numerical Approach}, leading to semi-discrete formulations for the advection-diffusion equation and the incompressible Navier-Stokes equations.
Section~\ref{sec:fully_discrete_formulations} presents the corresponding fully-discrete formulations and discusses various algorithmic aspects of the resulting numerical schemes. 
The performance of the method is assessed for various benchmark tests in Section~\ref{sec: Numerical Examples}. The tests demonstrate, for an appropriate choice of basis functions, optimal second and third-order convergence rates in space, and, given the particular time stepping strategy, second-order accuracy is obtained in time. Conclusions and an outlook for future research are presented in Section~\ref{sec: Conclusions Outlook}. 
\section{Problem formulation and definitions} \label{sec:problem-definitions}
We define the scalar-valued linear advection-diffusion equation and the incompressible Navier-Stokes equations on a domain $\Omega \subset \mathbb{R}^d$ (with $d=2,3$), having Lipschitz continuous boundary $\varGamma = \partial \Omega$. The boundary $\varGamma$ is partitioned into Dirichlet and  Neumann boundaries $\varGamma_D$  and $\varGamma_N$, respectively, satisfying $\varGamma_D \cup \varGamma_N = \partial\Omega$ and $\varGamma_D \cap \varGamma_N  = \emptyset$. The outward pointing unit vector normal to $\varGamma$ is denoted by $\mathbf{n}$. The time interval of interest is $I=\left(t^0, t^N \right]$, where $t^0$ and $t^N$ are the start and end time of the simulation, respectively. 
\subsection{Advection-diffusion equation}
On the space-time domain $\Omega \times I$, the scalar-valued linear advection-diffusion equation is defined as a system of two first-order equations as follows: given a solenoidal velocity field $\mathbf{a} \colon \Omega \times I \rightarrow \Real^d$, an initial condition $\phi_0: \Omega \rightarrow \Real$, diffusivity $\kappa$, and boundary conditions $h: \varGamma_N \times I \rightarrow \Real$ and $g: \varGamma_D \times I \rightarrow \Real$, find the scalar quantity $\phi:\Omega\times I \rightarrow \Real$ such that
	\begin{subequations} \label{eq:scalar-advection-diffusion_continuous}
		\begin{align}
		&\partialDt{\phi} + \Divergence{\sigma} = f  && \InDomain{\Omega \times I},  \label{eq:scalar_continuous}\\
		& \boldsymbol{\sigma}  = \mathbf{a}\phi - \kappa \nabla \phi && \InDomain{\Omega \times I},  \label{eq:flux_continuous}\\
		& \boldsymbol{\sigma} \cdot \mathbf{n} = \left(1-\gamma \right) (\mathbf{a} \cdot \mathbf{n})  \phi + h && \OnBoundary{\NeumannBC \times I}, \label{eq:scalar_neumann} \\ 
		&\phi = g  && \OnBoundary{\DirichletBC\times I} \label{eq:scalar_dirichlet}, \\
		& \phi(\mathbf{x},t^0) = \phi_0 && \InDomain{\Omega}. \label{eq:scalar_initial_condition}
		\end{align}
	\end{subequations}
The parameter $\gamma$ in Eq.~\eqref{eq:scalar_neumann} is equal to one at inflow boundaries (where $\mathbf{a} \cdot \mathbf{n} < 0$) and equal to zero on outflow Neumann boundaries (where $\mathbf{a}\cdot \mathbf{n} \geq 0$), with $h$ specifying the total flux on inflow Neumann boundaries and the diffusive flux on outflow Neumann boundaries.
\subsection{Incompressible Navier-Stokes equations}
The incompressible Navier-Stokes equations on the space-time domain $\Omega \times I$ are stated as follows: given the kinematic viscosity $\nu$, the forcing term $\mathbf{f}:\Omega\times I \rightarrow \mathbb{R}^d$, the boundary conditions $\mathbf{g}:\varGamma_D \times I \rightarrow \mathbb{R}^d$ and $\mathbf{h}:\varGamma_N \times I \rightarrow \mathbb{R}^d$, and a solenoidal initial condition $\boldsymbol{u}_0:\Omega \rightarrow \mathbb{R}^d$, find the velocity field $\boldsymbol{u}:\Omega \times I \rightarrow \mathbb{R}^d$ and pressure field $p:\Omega \times I \rightarrow \mathbb{R}$ such that
	\begin{subequations} \label{eq:incompressible-navier-stokes_continuous}
		\begin{align}
		&\partialDt{\boldsymbol{u}} +  \nabla \cdot \boldsymbol{\sigma} = \mathbf{f} &\text{in} \hspace{5pt}\Omega \times I,  \label{eq:navier-stokes_momentum}\\
		&\boldsymbol{\sigma}
		= 
		\boldsymbol{u}\otimes \boldsymbol{u} 
		+ 
		p\mathbf{I} 
		- 
		2\nu \nabla^s \boldsymbol{u} & \text{in} \hspace{5pt}\Omega \times I, \label{eq:navier-stokes_flux} \\
		&\nabla \cdot \boldsymbol{u} = 0 & \text{in} \hspace{5pt}\Omega \times I, \label{eq:navier-stokes_incompressibility} \\
		& \boldsymbol{\sigma} \mathbf{n}= \left(1-\gamma \right) \left(\boldsymbol{u} \otimes \boldsymbol{u} \right) \mathbf{n} + \mathbf{h} &\text{on} \hspace{5pt}\varGamma_N \times I, \label{eq:navier-stokes_neumann} \\
		& \boldsymbol{u} = \mathbf{g} & \text{on} \hspace{5pt} \varGamma_D \times I, \label{eq:navier-stokes_dirichlet} \\
		& \boldsymbol{u}(\mathbf{x},t^0) = \boldsymbol{u}_0(\mathbf{x}) & \text{in} \hspace{5pt} \Omega,  \label{eq:navier-stokes_initial-condition}
		\end{align} 
	\end{subequations}
where $\boldsymbol{\sigma}$ is the momentum flux tensor, $\mathbf{I}$ is the identity tensor, and $\nabla^s\boldsymbol{u}=\frac{1}{2}\nabla(\boldsymbol{u})+\frac{1}{2}\nabla(\boldsymbol{u})^\top$ is the symmetric velocity-gradient tensor. Similar to its use in the advection-diffusion equation, the parameter $\gamma$ is equal to one on inflow parts of $\varGamma_N$ (that is, where $\boldsymbol{u}\cdot \mathbf{n} < 0$), imposing in this way the total momentum flux. On the outflow parts of $\varGamma_N$ (where $\boldsymbol{u}\cdot \mathbf{n} \geq 0$) $\gamma$ equals zero, thus prescribing the diffusive part of the momentum flux only.
\subsection{Operator splitting} \label{eq:operator_splitting}
As in \cite{Maljaars2017}, the particle-mesh method is conceived as an operator splitting procedure. To this end, let the time interval of interest $I$ be partitioned using a sequence of $N+1$ discrete time levels $\{t^0, t^1,\dots,t^{N-1}, t^{N} \}$ which for $n=0,N-1$ define the half-open subintervals $I_n=(t^n,t^{n+1}]$ such that $\bigcup_n I_n=I$, while $\mathcal{I} := \left\{I_n\right\}$  defines the ordered sequence of subintervals. Furthermore, let the total fluxes (given by Eqs.~\eqref{eq:flux_continuous} and \eqref{eq:navier-stokes_flux}) be decomposed into advective parts $ \boldsymbol{\sigma}_a $ and diffusive parts $ \boldsymbol{\sigma}_d $.

A spatial-temporal operator splitting procedure for the advection-diffusion problem, Eq.~\eqref{eq:scalar-advection-diffusion_continuous}, now involves a scalar field $\psi: \Omega \times I_n \rightarrow \mathbb{R}$ satisfying an advection problem,
	\begin{subequations} \label{eq:scalar_hyperbolic_subproblem}
	\begin{align}
	&\partialDt{\psi} 
	+ \Divergence{\sigma}_a 
	= 0  && \InDomain{\Omega \times I_n},  \label{eq:scalar_hyperbolic_subproblem_transport}\\
	& \boldsymbol{\sigma}_a  = \mathbf{a}\psi && \InDomain{\Omega \times I_n},   \label{eq:scalar_hyperbolic_subproblem_flux}\\
	& \boldsymbol{\sigma}_a \cdot \mathbf{n} = (1-\gamma)(\mathbf{a} \cdot \mathbf{n}) \psi + \gamma h_a && \OnBoundary{\NeumannBC \times I_n}, 
	\label{eq:scalar_hyperbolic_subproblem_neumann} \\
	&\psi = g  && \OnBoundary{\DirichletBC^- \times I_n}, \label{eq:scalar_hyperbolic_subproblem_dirichlet}\\
	& \psi(\mathbf{x},t^n) = \mathcal{P}_L\left(\mathbf{\phi}(\mathbf{x},t^n) \right)&& \InDomain{\Omega},  \label{eq:scalar_mesh-particle_abstract}
	\end{align}
	\end{subequations}
and a scalar field $\phi: \Omega \times I_n \rightarrow \mathbb{R}$ satisfying a diffusion problem,
	\begin{subequations} \label{eq:diffusion_subproblem}
	\begin{align}
	&\partialDt{\phi} + \nabla \cdot \boldsymbol{\sigma}_d = f &&\InDomain{\Omega \times I_n}, \label{eq:scalar_diffusion_subproblem_continuity} \\ 
	&\boldsymbol{\sigma}_d = -\kappa \nabla \phi  && \InDomain{\Omega \times I_n}, \label{eq:scalar_diffusion_subproblem_flux}\\
	&\boldsymbol{\sigma}_d \cdot \mathbf{n} = h_d  && \OnBoundary{\NeumannBC \times I_n}, \\
	&\phi = g 									 && \OnBoundary{\DirichletBC \times I_n}, \label{eq:scalar_diffusion_subproblem_dirichlet} \\
	& \phi(\mathbf{x},t^{n}) = \mathcal{P}_E\left( \mathbf{\psi}(\mathbf{x},t^{n+1})\right) && \InDomain{\Omega}. \label{eq:scalar_particle-mesh_abstract} 
	\end{align}
	\end{subequations}
to be applied sequentially for every $I_n \in \mathcal{I}$.
In the advection stage the Dirichlet boundary condition can only be prescribed on inflow Dirichlet boundaries, denoted with $\varGamma_D^-$.
The flux prescribed on $\varGamma_N$ is split into an advective flux $h_a$ and a diffusive flux $h_d$, such that $h = h_a + h_d$. Note that the advective flux cannot be specified at outflow Neumann boundaries, which is automatically taken care of by virtue of Eq.~\eqref{eq:scalar_hyperbolic_subproblem_neumann}. 
Furthermore,  $\mathcal{P}_L$ and $\mathcal{P}_E$ are projection operators, which are introduced in order to couple the fields $\psi$ and $\phi$, with these fields being naturally defined on the particles and the mesh, respectively, in PIC methods.
More precisely, the projection operator $\mathcal{P}_L$ provides the initial condition at $t^n$ to advance the Lagrangian advection problem to $t^{n+1}$, and the projection operator $\mathcal{P}_E$ provides the initial condition at $t^n$ to advance the Eulerian diffusion problem to $t^{n+1}$. To remain consistent with Eq.~\eqref{eq:scalar-advection-diffusion_continuous} it is required that $\mathcal{P}_E \circ \mathcal{P}_L$ equals the identity operator. 

In a similar fashion, an operator splitting procedure for the incompressible Navier-Stokes equations (Eq.~\eqref{eq:incompressible-navier-stokes_continuous}) is formulated by finding a vector-valued field $\ContVector{v}: \Omega \times I_n \rightarrow \mathbb{R}^d$ satisfying an advection problem,
\begin{subequations} \label{eq:navier-stokes_advection_subproblem}
	\begin{align}
	&\partialDt{\ContVector{v}} +  \nabla \cdot \ContVector{\sigma}_a = \mathbf{0} & \text{in} \hspace{5pt}\Omega \times I_n, \label{eq:navier-stokes_advection_transport} \\
	&\ContVector{\sigma}_a = \ContVector{v} \otimes \ContVector{U} & \text{in} \hspace{5pt}\Omega \times I_n, \label{eq:navier-stokes_advection_subproblem_flux} \\
	&\ContVector{\sigma}_a \mathbf{n} = (1-\gamma)\left( \ContVector{v} \otimes \ContVector{U}\right) \mathbf{n} +  \gamma \mathbf{h}_a &\text{on} \hspace{5pt}\varGamma_N \times I_n, \label{eq:navier-stokes_advection_subproblem_neumann} \\
	& \boldsymbol{v} = \mathbf{g} & \text{on} \hspace{5pt} \varGamma_D^{-} \times I_n, \label{eq:navier-stokes_advection_subproblem_dirichlet} \\
	& \boldsymbol{v}(\mathbf{x},t^n) = \mathcal{P}_L\left(\boldsymbol{u}(\mathbf{x},t^n) \right)& \text{in} \hspace{5pt} \Omega,  \label{eq:vector_mesh-particle_abstract}
	\end{align}
\end{subequations}
and a velocity field $\boldsymbol{u}: \Omega \times I_n \rightarrow \mathbb{R}^d$ satisfying an
incompressible Stokes problem,
\begin{subequations} \label{eq:navier-stokes_stokes_subproblem}
	\begin{align} 
	&\partialDt{\boldsymbol{u}} +  \nabla \cdot \boldsymbol{\sigma}_d = \mathbf{f} &\text{in} \hspace{5pt}\Omega \times I_n,  \label{eq:navier-stokes_stokes_subproblem_momentum}\\
	&\nabla \cdot \boldsymbol{u} = 0 & \text{in} \hspace{5pt}\Omega \times I_n, \label{eq:navier-stokes_stokes_subproblem_incompressibility} \\
	&\boldsymbol{\sigma}_d = p\mathbf{I} - 2\nu \nabla^s \boldsymbol{u} & \text{in} \hspace{5pt}\Omega \times I_n, \label{eq:navier-stokes_stokes_subproblem_flux} \\
	& \boldsymbol{u} = \mathbf{g} & \text{on} \hspace{5pt} \varGamma_D \times I_n, \label{eq:navier-stokes_stokes_subproblem_dirichlet} \\
	& \boldsymbol{\sigma}_d\mathbf{n} = \mathbf{h}_d &\text{on} \hspace{5pt}\varGamma_N \times I_n, \label{eq:navier-stokes_stokes_subproblem_neumann} \\
	& \boldsymbol{u}(\mathbf{x},t^n) = \mathcal{P}_E\left( \boldsymbol{v}(\mathbf{x},t^{n+1})\right) & \text{in} \hspace{5pt} \Omega, \label{eq:vector_particle-mesh_abstract}
	\end{align}
\end{subequations}
with notation similar to Eqs.~(\ref{eq:scalar_hyperbolic_subproblem}-\ref{eq:diffusion_subproblem}). The advective field $\ContVector{U}$ in Eq.~\eqref{eq:navier-stokes_advection_subproblem_flux} is not yet specified, other than to require this field to be a consistent approximation to $\ContVector{u}$ which is piecewise constant on $\mathcal{I}$. 

We aim to formulate the projection operators $\mathcal{P}_L$ (in Eqs.~\eqref{eq:scalar_mesh-particle_abstract} and~\eqref{eq:vector_mesh-particle_abstract}) and $\mathcal{P}_E$ (in Eqs.~\eqref{eq:scalar_particle-mesh_abstract} and~\eqref{eq:vector_particle-mesh_abstract}) such that exact conservation is guaranteed, in the sense that the projected field satisfies Eq.~\eqref{eq:scalar_hyperbolic_subproblem_transport} or Eq.~\eqref{eq:navier-stokes_advection_transport} in an integral sense over each cell. 
\subsection{Definitions}
We next introduce some notation needed for the spatial discretizations of the Eulerian and Lagrangian subproblems, respectively, and the definition of the associated projections $\mathcal{P}_E$ and $\mathcal{P}_L$.
The fixed Eulerian mesh consists of the triangulation $\mathcal{T}:= \{K\}$ of the domain $\Omega$ into open, non-overlapping cells $K$. A measure of the cell size is denoted by $h_K$, and the outward pointing unit normal vector on the boundary $\partial K$ of each cell is denoted by $\mathbf{n}$. The closure of a cell is denoted by $\overline{K} = K \cup \partial K$. Adjacent cells $K_i$ and $K_j$ ($i\neq j$) share a common facet 
$F =\partial K_i \cap \partial K_j$. The set of all facets 
(including the exterior boundary facets $F =\partial K \cap \partial \Omega$) is denoted by $\mathcal{F}$. 
\subsubsection{Function spaces}
The following scalar finite element spaces are defined:
	\begin{align}
    W_h   &:= \left\{w_h \in L^2(\mathcal{T}), \hspace{3pt} w_h\lvert_K \hspace{3pt}  \in P_k(K) \hspace{3pt} \forall \hspace{3pt} K \in \mathcal{T} \right\}, \label{eq:wspace_local} \\
    T_h   &:= \left\{\tau_h \in L^2(\mathcal{T}), \hspace{3pt} \tau_h \lvert_K \hspace{3pt}  \in P_l(K) \hspace{3pt} \forall \hspace{3pt} K \in \mathcal{T} \right\}, \label{eq:lagrange_space} \\
	\bar{W}_{h,g} &:= \left\{\bar{w}_h \in L^{2}(\mathcal{F}), \hspace{3pt} \bar{w}_h \lvert_F \hspace{3pt} \in P_{k} (F) \hspace{3pt} \forall \hspace{3pt} F \in \mathcal{F}, \hspace{3pt} \bar{w}_h = g \hspace{3pt} \text{on} \hspace{3pt} \varGamma_D \right\},  \label{eq:wspace_global} \\
	Q_h &:= \left\{q_h \in L^2(\mathcal{T}), \hspace{3pt} q_h\lvert_K \hspace{3pt} \in P_{k-1}(K) \hspace{3pt} \forall \hspace{3pt} K \in \mathcal{T} \right\}, \label{eq:qspace_local}\\
	\bar{Q}_h &:= \left\{ \bar{q}_h \in L^{2}(\mathcal{F}), \hspace{3pt} \bar{q}_h \lvert_F \hspace{3pt} \in P_{k} (F) \hspace{3pt} \forall \hspace{3pt} F \in \mathcal{F}\right\},  \label{eq:qspace_global}
	\end{align}
in which $P(K)$ and $P(F)$ denote the spaces spanned by Lagrange polynomials on $K$ and $F$, respectively, and $k\geq 1 $ and $l\geq 0$ indicate the polynomial orders. Note that $\bar{W}_{h,g}$ and $\bar{Q}_{h}$ are defined only on cell facets and 
consist of functions which are continuous (i.e. single-valued) inside the facets $F \in \mathcal{F}$ and discontinuous at their borders. Furthermore, the facet function space $\bar{W}_{h,g}$ satisfies the inhomogeneous Dirichlet boundary condition on $\varGamma_D$, with the related space $\bar{W}_{h,0}$ satisfying the homogeneous Dirichlet boundary condition on $\Gamma_D$. 
Importantly, the extension of $W_h, T_h, Q_h$ to the cell boundary $\partial K$ is formally only defined via a trace operator, however, we omit this technicality in the sequel to avoid notational clutter. 

Finally, $\boldsymbol{W}_h$, $\boldsymbol{T}_h$, and $\boldsymbol{\bar{W}}_{h,g}$ denote the finite element spaces of $d$-vectors in $\mathbb{R}^d$ corresponding to the scalar function spaces $W_h$,  $T_h$, and $\bar{W}_{h,g}$, respectively.
\subsubsection{Particle definitions}
The Lagrangian particle configuration in the domain $ \Omega $ at a fixed time instant $ t $ is defined as follows
	\begin{equation}
	\mathcal{X}_t := \{\mathbf{x}_p(t) \in \Omega\}_{p=1}^{N_p},
	\end{equation}
in which $\mathbf{x}_p$ denotes the spatial coordinates of particle $ p $ and $ N_p $ is the number of particles. \\
Furthermore, the index set of all particles and the index set of particles hosted by cell $K$, at a fixed time instant $t$ are defined as, respectively 
	\begin{align} 
	\mathcal{S}_t &:= \defineSet{p\in\mathbb{N}}{\mathbf{x}_p(t) \in \mathcal{X}_t}, \label{eq:particle set definition} \\
	\mathcal{S}^{K}_t &:= \defineSet{p\in \mathbb{N}}{\mathbf{x}_p(t) \in \overline{K}, \hspace{3pt} \mathbf{x}_p(t) \in \mathcal{X}_t}. \label{eq:particle set definition_local}
	\end{align} 
Finally, Lagrangian scalar and vector fields on the particles are defined as, respectively,
	\begin{align}
	\Psi_t &:= \left\{ \psi_p(t) \in \mathbb{R}, \hspace{3pt} \forall \hspace{3pt} p \in \mathcal{S}_t  \right \}, \\
	\mathcal{V}_t &:= \left\{ \mathbf{v}_p(t) \in \mathbb{R}^d, \hspace{3pt} \forall \hspace{3pt} p \in \mathcal{S}_t \right \},
	\end{align}
where $\psi_p$ and $\mathbf{v}_p(t)$ denote the corresponding scalar and vector quantities associated with particle $p$.\\
Importantly, subscripts $p$ and $h$ are consistently used to distinguish between Lagrangian particle data and Eulerian mesh fields, respectively. 
\section{Semi-discrete formulations} \label{sec: Numerical Approach}
We provide now a collection of semi-discrete forms that together constitute a particle-mesh operator splitting method for the advection-diffusion equation (Eqs.~(\ref{eq:scalar_hyperbolic_subproblem}-\ref{eq:diffusion_subproblem})) and for the incompressible Navier-Stokes equations (Eqs.~(\ref{eq:navier-stokes_advection_subproblem}-\ref{eq:navier-stokes_stokes_subproblem})). 
The splitting procedure involves the following sequence of steps (referring to the relevant equations for the advection-diffusion operator splitting for the sake of illustration):
\begin{enumerate}[noitemsep]
\item \textit{Lagrangian discretization of the advection problem}, in order to solve Eqs.~(\ref{eq:scalar_hyperbolic_subproblem_transport}~-~\ref{eq:scalar_hyperbolic_subproblem_dirichlet}) at the particle level;
\item \textit{particle-mesh projection}, in order to project the scattered particle data onto a (scalar) field on the Eulerian mesh using the operator $\mathcal{P}_E: \Psi_t \rightarrow W_h$  (see Eq.~\eqref{eq:scalar_particle-mesh_abstract});
\item \textit{Eulerian discretization of the diffusion equation}, in order to solve Eqs.~(\ref{eq:scalar_diffusion_subproblem_continuity}~-~\ref{eq:scalar_diffusion_subproblem_dirichlet}) on the mesh;
\item \textit{mesh-particle projection}, in order to update the particle properties 
from the scalar field on the Eulerian mesh
using the operator $\mathcal{P}_L: W_h \rightarrow \Psi_t$  (see Eq.~\eqref{eq:scalar_mesh-particle_abstract}). 
\end{enumerate} 

The resulting semi-discrete components for the advection-diffusion problem  and the incompressible Navier-Stokes problem are presented in Sections~\ref{sec:adv-diffusion_semi_discrete} and \ref{sec:navier-stokes_semi-discrete}, respectively. Since the formulation of conservative 
projection operators is our main objective we primarily focus on steps 2 and 4, referring to Maljaars~et~al.~\cite{Maljaars2017} for a discussion of steps 1 and 3. 
\subsection{Advection-diffusion equation} \label{sec:adv-diffusion_semi_discrete}
\subsubsection{Lagrangian discretization of the advection problem}
In a Lagrangian, particle-based frame of reference, solving the advection problem Eqs.~(\ref{eq:scalar_hyperbolic_subproblem_transport}-\ref{eq:scalar_hyperbolic_subproblem_dirichlet}) on the time interval $I_n$ is straightforward by decomposing the problem into two ordinary differential equations for the particle scalar quantity and the particle position, given by
	\begin{subequations} \label{eq:lagrangian_advection_linear}
	\begin{align}
	\dot{\psi}_p(t) &= 0 && \forall \hspace{3pt} p \in \mathcal{S}_t , \label{eq5:particle_vel_update} \\
	\dot{\mathbf{x}}_p(t)&=  
	\mathbf{a}(\mathbf{x}_p(t),t^n) && \forall \hspace{3pt} p \in \mathcal{S}_t,  \label{eq5:particle_advection}
	\end{align}
	\end{subequations}
where $\dot{\psi}_p(t)$ and $\dot{\mathbf{x}}_p(t)$ are the total derivatives at time $t \in I_n$ of 
the scalar quantity and the position of particle $p$, respectively. Furthermore, $\mathbf{a}(\mathbf{x},t^n)$ is a prescribed solenoidal velocity field at time $t^n$. An important observation is that $\psi_p$ stays constant throughout the particle advection stage by virtue of Eq.~\eqref{eq5:particle_vel_update}.
\subsubsection{PDE-constrained particle-mesh interaction} \label{sec:pde_constrained_interaction}
In \cite{Maljaars2017} the projection operator $\mathcal{P}_E$, transferring the scattered particle data to piecewise continuous fields at the mesh, is formulated in terms of local least-squares projections obtained from the minimization problem
\begin{equation}
    \min_{\psi_h \in W_{h}} \, J = \sum_{p \in \mathcal{S}_t}^{} \frac{1}{2}\left( \psi_h(\mathbf{x}_p(t), t) - \psi_p(t)\right)^2. 
    \label{eq:J1}
\end{equation}
With $W_h$ a discontinuous function space, this approach allows for an efficient cellwise implementation, and gives accurate results provided that the particle locations satisfy unisolvency with respect to $W_h$. However,
conservation of linear quantities is not guaranteed after projection. We now extend the variational framework of that approach by imposing constraints as to ensure conservation throughout the particle-mesh projection.
Recalling that the Lagrangian particles are used for solving an advection operator, it is argued that the projected particle motion should also satisfy an advection operator from the perspective of the Eulerian mesh.
Therefore, we constrain the particle-mesh projection by a hyperbolic conservation law. 

To this end, we augment the functional in Eq.~\eqref{eq:J1} with terms multiplying Eq.~\eqref{eq:scalar_hyperbolic_subproblem} with a Lagrange multiplier $\lambda_h \in T_h$. Integration by parts leaves an unknown flux on interior cell facets which is formulated in terms of a control variable $\Bar{\psi}_h \in \Bar{W}_{h,g}$, while the Neumann boundary condition is substituted on exterior facets. For a given particle field $\psi_p \in \Psi_t$, an advective velocity field $\mathbf{a} : \Omega \times I_n \rightarrow \mathbb{R}^d$, the initial condition $\psi_h^n \in W_h $, and an advective Neumann boundary condition $h_a: \Gamma_N \times I_n \rightarrow \mathbb{R}$, the minimization problem then involves finding the stationary points of the Lagrangian functional 
	\begin{multline} \label{eq:scalar_lagrangian-functional}
	\mathcal{L}(\psi_h, \bar{\psi}_h, \lambda_h) 
	= 
	\sum_{p \in \mathcal{S}_t}^{} \frac{1}{2}\left( \psi_h(\mathbf{x}_p(t), t) - \psi_p(t)\right)^2 
	+ 
	\sum\limits_{K}^{}\lineIntegral{\partial K}{\frac{1}{2}\beta \left( \bar{\psi}_h - \psi_h \right)^2} 
	+ 
	\areaIntegral{\Omega}{\partialDt{\psi_h}\lambda_h }  \\
	- 
	\sum\limits_{K}^{} \areaIntegral{K}{ \mathbf{a} \psi_h \cdot \GradientScalar{\lambda_h} } 
	+ 
	\sum\limits_{K}^{} \lineIntegral{\partial K \setminus \NeumannBC}{ \mathbf{a}\cdot \mathbf{n} \bar{\psi}_h \lambda_h} 
	+
	\lineIntegral{\NeumannBC}{ \left(1-\gamma\right) \mathbf{a}\cdot \mathbf{n} \xDiscreteScalar{\psi} \xDiscreteScalar{\lambda}} 
	+ 
	\lineIntegral{\NeumannBC}{\gamma h_a \lambda_h },  
	\end{multline}
for every $ t \in I_n$. The collection of terms containing $\lambda_h$ constitutes a weak form of the advection subproblem, Eq.~\eqref{eq:scalar_hyperbolic_subproblem}. Furthermore, the unknown facet-based field $\bar{\psi}_h \in \bar{W}_{h,g}$ determines the interface flux, the implications of which will become clear later on. The additional term containing $\beta > 0$ penalizes the jumps between $\psi_h$ and $\Bar{\psi}_h$ on cell interfaces thereby avoiding the problem of becoming singular in cases with vanishing normal velocity $\mathbf{a} \cdot \mathbf{n}$.

Equating the variations of Eq.~\eqref{eq:scalar_lagrangian-functional} with respect to the three unknowns $\left(\psi_h, \lambda_h, \bar{\psi}_h \right) \in \left(W_h, T_h, \bar{W}_{h,g} \right)$ to zero, results in the following system of variational equations. At time $t \in I_n$, variation with respect to the scalar field $\xDiscreteScalar{\psi}$ gives the co-state equation
\begin{subequations}\label{eq:scalar_optimality_conditions}
    \begin{align}
		\begin{split}
		{}& \sum_{p \in \mathcal{S}_t} \left( \psi_h(\mathbf{x}_p(t), t) - \psi_p(t) \right) \delta \psi_h (\mathbf{x}_p(t))
		 -
		 \sum\limits_{K}^{} \lineIntegral{\partial K}{\beta \left( \bar{\psi}_h - \psi_h \right) \delta \psi_h}
		 \\
		& \hspace{15pt} 
		+ 
		\areaIntegral{\Omega}{ \frac{\partial \delta \psi_h}{\partial t} \lambda_h}  
		-   
		\sum_{K}\areaIntegral{K}{ \mathbf{a} \cdot \GradientScalar{\lambda_h} \delta \psi_h} 
		+
		\lineIntegral{\NeumannBC}{ \left(1-\gamma\right) \mathbf{a}\cdot \mathbf{n} \xDiscreteScalar{\lambda} \delta \xDiscreteScalar{\psi} } 
		= 
		0 \hspace{10pt } \forall \hspace{3pt} \delta \psi_h \in W_h.
		\end{split} \label{eq:objective_semi_discrete}
	\end{align}
	Variation with respect to the Lagrange multiplier $\xDiscreteScalar{\lambda}$ gives the state equation, 
	\begin{align}
	\begin{split}
		&\areaIntegral{\Omega}{\frac{\partial \psi_h }{\partial t} \delta \lambda_h }    - 
		 \sum_{K}\areaIntegral{K}{\mathbf{a} \psi_h   \cdot \GradientScalar{\delta \lambda_h} }  
		+
		\sum_{K}^{}\lineIntegral{\partial K\setminus \varGamma_N}{ \mathbf{a} \cdot \mathbf{n}  \bar{\psi}_h \delta \lambda_h } \\
		& \hspace{50pt}
		+
		\lineIntegral{\NeumannBC}{ \left(1-\gamma\right) \mathbf{a}\cdot \mathbf{n} \xDiscreteScalar{\psi} \delta \xDiscreteScalar{\lambda}} 
		+ 
		\lineIntegral{\varGamma_N}{ \gamma h_a \delta \lambda_h}  
		=
		0 \hspace{10pt} \forall \hspace{3pt} \delta \lambda_h \in T_h.
		\end{split}  \label{eq:semi_discr_state_eq}
	\end{align}
	Finally, variation with respect to the control variable $\xDiscreteScalar{\bar{\psi}}$ leads to the optimality condition,
	\begin{equation}
	\sum\limits_{K}^{} \lineIntegral{\partial K\setminus \varGamma_N}{ \mathbf{a}\cdot \mathbf{n} \lambda_h \delta \bar{\psi}_h}
		+
		\sum\limits_{K}^{} \lineIntegral{\partial K}{\beta \left( \bar{\psi}_h - \psi_h \right) \delta \bar{\psi}_h } = 0  \hspace{30pt}  \forall \hspace{3pt} \delta \bar{\psi}_h \in \bar{W}_{h,0}. \label{eq:optimality condition}
	\end{equation}
\end{subequations}

Note that the optimality system Eq.~\eqref{eq:scalar_optimality_conditions} is naturally formulated using a HDG approach, since we require the optimal solution for the state variable $\xDiscreteScalar{\psi}$ to be in the discontinuous function space $\xDiscreteScalar{W}$, and the optimal control to be provided by the single-valued, facet-based control variable $\xDiscreteScalar{\bar{\psi}} \in \bar{W}_{h,g}$.
After an appropriate discretization of the time derivatives in Eqs.~(\ref{eq:objective_semi_discrete} and~\ref{eq:semi_discr_state_eq}) - to be discussed in Section~\ref{sec:fully_discrete_formulations} - a field $\xDiscreteScalar{\psi} \in \xDiscreteScalar{W}$ can be reconstructed from the particle data $\psi_p \in \Psi_t$ by solving the optimality system Eq.~\eqref{eq:scalar_optimality_conditions}. This reconstructed field is used to provide the initial condition for the subsequent diffusion subproblem.
\subsubsection{Eulerian discretization of the diffusion equation}
We next present the semi-discrete formulation of the diffusion problem Eq.~\eqref{eq:diffusion_subproblem}. To this end, we will use the HDG method presented by Labeur~\&~Wells~\cite{Labeur2007} to seek solutions $\phi_h \in W_h$. This choice allows for trivial projections between fields $\psi_h$ and $\phi_h$, which seamlessly fits in the approach used for solving the optimality system Eq.~\eqref{eq:scalar_optimality_conditions}. 
Referring to \cite{Labeur2007} for further details, the HDG discretization results in a set of local and global problems which are respectively stated as: at time $t\in I_n$, given the initial condition $\mathbf{\phi}_h^{n} \in W_h$, the diffusive Neumann boundary condition $h_d : \varGamma_N \rightarrow \mathbb{R}$ and the diffusivity $\kappa$, find $ \xDiscreteScalar{\phi}\in W_h$ and  $\xDiscreteScalar{\bar{\phi}} \in \bar{W}_{h,g}$ such that
	\begin{subequations} \label{eq:diffusion_semi_discrete}
	\begin{equation}
	\begin{alignedat}{2}
	\areaIntegral{\Omega}{\partialDt{\xDiscreteScalar{\phi}} \xDiscreteScalar{w}} \;
	+ \;
	&\sum\limits_{K}^{}\areaIntegral{K}{\kappa \nabla \xDiscreteScalar{\phi}\cdot \nabla \xDiscreteScalar{w}}
	\; + \;
	\sum\limits_{K}^{}\lineIntegral{\partial K}{
	\boldsymbol{\Hat{\sigma}}_{d,h} \cdot \mathbf{n} \xDiscreteScalar{w}} 
	\\
	& \qquad +
	\sum\limits_{K}^{}\lineIntegral{\partial K}{\kappa \left(\xDiscreteScalar{\bar{\phi}} - \xDiscreteScalar{\phi} \right)\mathbf{n} \cdot \nabla \xDiscreteScalar{w}}
	= 
	\areaIntegral{\Omega}{f \xDiscreteScalar{w}} && \qquad \forall \hspace{3pt} \xDiscreteScalar{w} \in W_h, \\
	\end{alignedat}
	\end{equation}
	and  
	\begin{align}
	\sum\limits_{K}^{}\lineIntegral{\partial K}{\boldsymbol{\Hat{\sigma}}_{d,h}  \cdot \mathbf{n} \xDiscreteScalar{\bar{w}}} 
	\; = \; 
	\lineIntegral{\NeumannBC}{h_d \xDiscreteScalar{\bar{w}}} && \forall \hspace{3pt} \xDiscreteScalar{\bar{w}}\in \bar{W}_{h,0},
	\end{align}
	where $\boldsymbol{\Hat{\sigma}}_{d,h}$ is a diffusive flux at cell facets defined by 
    \begin{align}
    \boldsymbol{\Hat{\sigma}}_{d,h} = - \kappa \nabla \phi - \frac{\alpha}{h_K} \kappa \left(\Bar{\phi}_h - \phi_h\right)\mathbf{n},
    \end{align}
	\end{subequations}
in which $\alpha$ is a dimensionless parameter as is typical to interior penalty methods \cite{Arnold2002}.
Solving Eq.~\eqref{eq:diffusion_semi_discrete} yields $\xDiscreteScalar{\bar{\phi}}$ and  $\xDiscreteScalar{\phi}$, where the latter is to be used for updating the particle properties $\psi_p$ in the subsequent mesh-particle projection step.
\subsubsection{Mesh-particle projection}
The mesh-particle projection $\mathcal{P}_L : W_h \rightarrow \Psi_t$ is based on the following minimization problem 
	\begin{equation} \label{eq:mesh2particle_semidiscrete}
	\min_{\psi_p(t)} \, J := \sum_{p\in \mathcal{S}_t}^{}\frac{1}{2} \left( \phi_h(\mathbf{x}_p(t),t) - \psi_p(t) \right)^2,
	\end{equation}
where we emphasize that the objective functional $J$ is also at the basis of the particle-mesh projection Eq.~\eqref{eq:scalar_lagrangian-functional}.
Carrying out the minimization yields the particularly simple result
	\begin{align} \label{eq:mesh-particle_semi-discrete}
	\psi_p(t) = \phi_h(\mathbf{x}_p(t),t)  && \forall \hspace{3pt} p \in \mathcal{S}_t.
	\end{align}
The mesh-particle projection is not restricted to the mapping of the scalar field $\phi_h \in W_h$ only, but can be applied to project arbitrary fields in $W_h$~-~e.g. the temporal increments of $\phi_h$, see Section~\ref{sec:mesh2particle_discrete}~-~onto the particles. 
As such, it provides the tool for updating the particle quantities, which completes the semi-discrete sequence of steps comprising the particle-mesh operator splitting of the advection-diffusion equation. 
\subsection{Incompressible Navier-Stokes equations} \label{sec:navier-stokes_semi-discrete}
In the same vein, we consider a particle-mesh operator splitting for the incompressible Navier-Stokes equations. Details of the semi-discrete operator splitting procedure, which involves a non-linear advection step followed by an unsteady Stokes step, can be found in \cite{Maljaars2017}. The main difference between the approach in the aforementioned reference and the method presented here, is the PDE-constrained particle-mesh projection performed in step 2. For completeness, we also briefly recapitulate the remaining model components. 

To set the stage, some additional notation is introduced. The momentum field at the Eulerian mesh, to be reconstructed from the particle field, is denoted by $\xDiscreteVector{v} \in \xDiscreteVector{W}$. Throughout the advection stage we will make use of advective fields $\left(\xDiscreteVector{U}, \xDiscreteVector{\Bar{U}} \right) \in \left( \xDiscreteVector{W}, \boldsymbol{\Bar{W}}_{h,g} \right)$.
Velocities and pressures at the Eulerian mesh related to the Stokes step are denoted by $(\xDiscreteVector{u}, \xDiscreteScalar{p}) \in \left(\xDiscreteVector{W}, \xDiscreteScalar{Q}\right)$. In addition, vector-valued facet functions $\xDiscreteVector{\bar{v}} \in  \boldsymbol{\bar{W}}_{h,g}$ and $\xDiscreteVector{\bar{u}} \in \boldsymbol{\bar{W}}_{h,g}$, and scalar-valued facet functions $\xDiscreteScalar{\bar{p}} \in \xDiscreteScalar{\bar{Q}}$ are used.
\subsubsection{Lagrangian discretization of the advection problem}
The particle advection step for the incompressible Navier-Stokes equations proceeds along similar lines as for the linear advection-diffusion equation (Eq.~\eqref{eq:lagrangian_advection_linear}). 
Particle specific momenta and positions are updated by integrating
    \begin{align} 
    \dot{\mathbf{v}}_p(t)
    &=
    \mathbf{0} && \forall \hspace{3pt} p \in \mathcal{S}_t, \label{eq5:ns_specific_momentum} \\
    \dot{\mathbf{x}}_p(t)
    &=  
    \xDiscreteVector{U}(\mathbf{x}_p(t),t^n) && \forall \hspace{3pt} p \in \mathcal{S}_t,  \label{eq5:ns_particle_advection}
    \end{align}
over the time interval $I_n$. In this equation $\mathbf{x}_p(t)$ denotes the particle position at time $t \in I_n$. We still leave the advective field $\xDiscreteVector{U} \in \xDiscreteVector{W}$ unspecified, but recall that it is a consistent approximation of $\mathbf{u}_h$ which is piecewise constant in $\mathcal{I}$. Similar to the scalar-valued case,  $\mathbf{v}_p$ stays constant throughout the advection stage by virtue of Eq.~\eqref{eq5:ns_specific_momentum}.
\subsubsection{PDE-constrained particle-mesh interaction}
The formulation of the vector-valued PDE-constrained projection operator $\mathcal{P}_E: \mathcal{V}_t \rightarrow \xDiscreteVector{W}$ of step 2 proceeds along similar lines as for the linear advection-diffusion equation. That is, the functional constituting the least-squares minimization problem (as in Eq.~\eqref{eq:J1}) is augmented with Lagrange multiplier terms that enforce weak satisfaction of the momentum advection equation by the projected field $\xDiscreteVector{v}$. Integration by parts of the constraint introduces a facet-based variable $\xDiscreteVector{\Bar{v}}$ controlling the flux across cell interfaces. The resulting optimality system can be stated as follows: at time $t\in I_n$, given the initial condition $\xDiscreteVector{v}^{n} \in \mathbf{W}_h$, the particle field $\boldsymbol{v}_p(t) \in \mathcal{V}_t$, the advective field $\left(\xDiscreteVector{U},\xDiscreteVector{\Bar{U}}\right) \in \left(\xDiscreteVector{W},\xDiscreteVector{\Bar{W}}\right)$ and the advective Neumann boundary condition $\mathbf{h}_a : \varGamma_N \rightarrow \mathbb{R}^d$, find $\left(\xDiscreteVector{v}, \xDiscreteVector{\lambda}, \xDiscreteVector{\bar{v}} \right) \in \left(\xDiscreteVector{W}, \xDiscreteVector{T}, \boldsymbol{\bar{W}}_{h,g}\right)$ such that 
\begin{subequations} \label{eq:optimality_NS}
	\begin{multline}\label{eq:costate-semi-discrete_NS}
	\sum_{p \in \mathcal{S}_t} \left( \xDiscreteVector{v}(\mathbf{x}_p(t), t) - \boldsymbol{v}_p \right) \cdot \delta \xDiscreteVector{v}(\mathbf{x}_p(t))   
	- 
	\sum_{K} \lineIntegral{\partial K}{\beta \left( \xDiscreteVector{\bar{v}} - \xDiscreteVector{v} \right) \cdot \delta \xDiscreteVector{v}} 
	+ 
	\areaIntegral{\Omega}{\partialDt{\delta \xDiscreteVector{v}} \cdot \xDiscreteVector{\lambda} }  \\  
	- 
	\sum_{K}\areaIntegral{K}{ \left(\delta \xDiscreteVector{v} \otimes \xDiscreteVector{U} \right) \colon \nabla \xDiscreteVector{\lambda} }  
	+
	\lineIntegral{\NeumannBC}{(1-\gamma) \left(\delta \xDiscreteVector{v} \otimes \xDiscreteVector{U}\right) \mathbf{n} \cdot \xDiscreteVector{\lambda}}
	= 
	0 \hspace{10pt} \forall \hspace{3pt}  \delta \xDiscreteVector{v} \in \xDiscreteVector{W},
	\end{multline}
	\begin{multline} \label{eq:state-semi-discrete_NS}
	\areaIntegral{\Omega}{\partialDt{\xDiscreteVector{v}} \cdot \delta \xDiscreteVector{\lambda}  }    
	- 
	\sum_{K}\areaIntegral{K}{\left(\xDiscreteVector{v} \otimes  \xDiscreteVector{U} \right)  : \nabla \delta \xDiscreteVector{\lambda}}  
	+
	\sum_{K}^{}\lineIntegral{\partial K\setminus \NeumannBC}{ \left( \xDiscreteVector{\bar{v}} \otimes \xDiscreteVector{\bar{U}} \right) \mathbf{n} \cdot \delta \xDiscreteVector{\lambda}}  
	\\
	+
	\lineIntegral{\NeumannBC}{(1-\gamma) \left(\xDiscreteVector{ v} \otimes \xDiscreteVector{U}\right) \mathbf{n} \cdot \delta \xDiscreteVector{\lambda}}
	+
	\oint\limits_{\NeumannBC}{ \gamma \mathbf{h}_a \cdot \delta \xDiscreteVector{\lambda}}  
	=
	0
	\hspace{10pt} \forall \hspace{3pt}  \delta \xDiscreteVector{\lambda} \in \xDiscreteVector{T}, 
	\end{multline}
	\begin{align}\label{eq:optimality-semi-discrete_NS}
	\sum_{K}^{}\lineIntegral{\partial K\setminus \NeumannBC}{\left( \delta \xDiscreteVector{\bar{v}} \otimes \xDiscreteVector{\bar{U}} \right) \mathbf{n} \cdot  \xDiscreteVector{\lambda} } 
	%
	+
	\sum\limits_{K}^{} \lineIntegral{\partial K}{\beta \left( \xDiscreteVector{\bar{v}} - \xDiscreteVector{v} \right) \cdot \delta \xDiscreteVector{\bar{v}} }
	= 0 
	&& \forall \hspace{3pt} \delta \xDiscreteVector{\bar{v}} \in \boldsymbol{\bar{W}}_{h,0}.
	\end{align}
\end{subequations}
It is worth mentioning the similarity of Eq.~\eqref{eq:optimality_NS} with Eq.~\eqref{eq:scalar_optimality_conditions} for the linear scalar advection problem, owing to the use of the advective fields $\xDiscreteVector{U}$ and $\xDiscreteVector{\bar{U}}$ which linearizes the momentum fluxes, thereby simplifying the optimality system considerably. \\
The PDE-constrained projection step leads to a field $\xDiscreteVector{v} \in \xDiscreteVector{W}$ to be used in the subsequent unsteady Stokes problem.
\subsubsection{Unsteady Stokes equations}
 The semi-discrete Stokes problem is formulated in terms of a velocity field 
$\xDiscreteVector{u} \in \xDiscreteVector{W}$ which enables a trivial projection between the fields $\xDiscreteVector{u}$ and $\xDiscreteVector{v}$. This practically implies that we will resort to the class of (H)DG methods for the spatial discretization of the unsteady Stokes equations. More specifically, the HDG method from Rhebergen~\&~Wells~\cite{Rhebergen2016, Rhebergen2017} is employed. In a particle-mesh setting, the restriction on the choice of function spaces advocated by these authors has some distinct advantages compared to the general formulation of Labeur~\&~Wells in earlier work \cite{Labeur2012}. This will be further discussed in Section~\ref{sec:function_space_general_remarks}.\\
Following \cite{Rhebergen2017}, the semi-discrete form of the unsteady incompressible Stokes problem on the time interval $I_n$ is stated as follows: at time $t \in I_n$, given the initial condition $\xDiscreteVector{u}^{n} \in \mathbf{W}_h$, a forcing term $ \ContVector{f} : \Omega \rightarrow \mathbb{R}^d$, the diffusive Neumann boundary condition  $\ContVector{h}_d : \varGamma_N \rightarrow \mathbb{R}^d$ and the viscosity $\nu$, find $\left( \xDiscreteVector{u},\xDiscreteVector{\bar{u}}, \xDiscreteScalar{p}, \xDiscreteScalar{\bar{p}} \right) \in \left( \xDiscreteVector{W}, \boldsymbol{\bar{W}}_{h,g}, Q_h, \bar{Q}_h \right)$ satisfying the local and global momentum balances, respectively,
	\begin{subequations} \label{eq:unsteady-stokes_semi-discrete}
	\begin{multline} \label{eq:sd local momentum}
	\areaIntegral{\Omega}{\partialDt{\xDiscreteVector{u}}\cdot \xDiscreteVector{w}}-\sum\limits_{K}^{}\areaIntegral{K}{\boldsymbol{\sigma}_{d,h} : \nabla \xDiscreteVector{w}}
	+ 
	\sum\limits_{K}^{}\lineIntegral{\partial K}{\boldsymbol{\hat{\sigma}}_{d,h}\mathbf{n}\cdot \xDiscreteVector{w}} 
	\\ 
	+
	\sum\limits_{K}^{}\lineIntegral{\partial K}{2\nu \left(\bar{\mathbf{u}}_h-\xDiscreteVector{u}\right)\cdot \nabla^s\xDiscreteVector{w}\mathbf{n}}
	= 
	\areaIntegral{\Omega}{\mathbf{f}\cdot \xDiscreteVector{w}} \hspace{15pt}\forall \hspace{3pt} \xDiscreteVector{w} \in \xDiscreteVector{W},
	\end{multline}
	\begin{align}\label{eq:sd global momentum}
	\sum\limits_{K}^{}\lineIntegral{\partial K}{\boldsymbol{\hat{\sigma}}_{d,h}\mathbf{n} \cdot \xDiscreteVector{\bar{w}}} 
	= 
	\oint\limits_{\varGamma_N}\mathbf{h}_d\cdot \xDiscreteVector{\bar{w}} \text{d}\varGamma && \forall \hspace{3pt} \xDiscreteVector{\bar{w}} \in \boldsymbol{\bar{W}}_{h,0},
	\end{align}
	and, simultaneously, the local and global mass balances, respectively, 
	\begin{align} \label{eq:sd local mass}
	 \sum\limits_{K}^{}\areaIntegral{K}{ \xDiscreteVector{u}\cdot \nabla q_h} 
	 - 
	 \sum\limits_{K}^{}\lineIntegral{\partial K}{\xDiscreteVector{u}\cdot \mathbf{n} \hspace{1.5pt} q_h} 
	 &= 
	 0 && \forall \hspace{3pt} q_h \in Q_h, \\
	 \sum\limits_{K}^{}\lineIntegral{\partial K}{\xDiscreteVector{u}\cdot\mathbf{n}\hspace{1.5pt} \bar{q}_h} 
	 - \lineIntegral{\partial \Omega}{\xDiscreteVector{\bar{u}}\cdot\mathbf{n} \hspace{1.5pt} \bar{q}_h} 
	 &= 
	 0  && \forall \hspace{3pt} \bar{q}_h \in \bar{Q}_h,
	\end{align}
	in which the diffusive interface flux $\boldsymbol{\hat{\sigma}}_{d,h}$ in Eq.~\eqref{eq:sd local momentum} and Eq.~\eqref{eq:sd global momentum} is defined as 
	\begin{equation} \label{eq:discrete momentum facet}
	    \boldsymbol{\hat{\sigma}}_{d,h} = \bar{p}_h \ContVector{I} - 2\nu\nabla^s \xDiscreteVector{u} - \frac{\alpha}{h_K}2\mathbf{\nu}\left(\xDiscreteVector{\bar{u}} - \xDiscreteVector{u}\right) \otimes \ContVector{n},
	\end{equation}
	\end{subequations}
with $\ContVector{I}$ the identity tensor, and $\alpha$ a parameter as is typical to interior penalty methods \cite{Arnold2002}. 

Solving Eq.~\eqref{eq:unsteady-stokes_semi-discrete} for $\left(\xDiscreteVector{u},\xDiscreteVector{\bar{u}}, \xDiscreteScalar{p},\xDiscreteScalar{\bar{p}} \right)$ yields the divergence-free velocity field $\xDiscreteVector{u} \in \xDiscreteVector{W}$ which is to be used in the subsequent mesh-particle projection step.
\subsubsection{Mesh-particle projection}
Next, the vector-valued counterpart of Eq.~\eqref{eq:mesh2particle_semidiscrete} is used for the mesh-particle projection $\mathcal{P}_L \colon \xDiscreteVector{W} \rightarrow \mathcal{V}_t$. Thus, the specific momentum $\mathbf{v}_p$ associated with particle $p$ is updated via
	\begin{align}
		\mathbf{v}_p(t) = \xDiscreteVector{u}\left(\mathbf{x}_p(t),t\right) && \forall p \in \mathcal{S}_t.
	\end{align}
This mesh-particle projection concludes the sequence of steps for a particle-mesh operator splitting of the incompressible Navier-Stokes equations.
\subsection{Properties of the semi-discrete formulation} \label{sec:properties_semi_discrete}
To better illustrate the properties of the PDE-constrained particle-mesh projection, we now prove consistency and conservation of the above presented scheme. 
In addition, we prove consistency and conservation of the resulting particle-mesh operator splitting schemes for the advection-diffusion equation and the incompressible Navier-Stokes equations. 
\subsubsection{Consistency} \label{sec:consistency1}
Consistency of the particle-mesh operator splitting scheme entails three different aspects, 
(i) the operator splitting has to be consistent with the unsplit governing equations (Eqs.~\eqref{eq:scalar-advection-diffusion_continuous} or \eqref{eq:incompressible-navier-stokes_continuous}), (ii) the Eulerian part (step 3) has to be consistent with the diffusion subproblem (Eqs.~\eqref{eq:diffusion_subproblem} or \eqref{eq:navier-stokes_stokes_subproblem}), and (iii), the constraint imposed weakly in the projection operator $\mathcal{P}_E$ (step 2) has to be consistent with the advection subproblem
(Eqs.~\eqref{eq:scalar_hyperbolic_subproblem} or \eqref{eq:navier-stokes_advection_subproblem}).

To start with (i), splitting of the advection-diffusion equation (Eq.~\eqref{eq:scalar-advection-diffusion_continuous}) and the incompressible Navier-Stokes equations (Eq.~\eqref{eq:incompressible-navier-stokes_continuous}) into a kinematic part (advection problem) and a constitutive part (diffusion equation or Stokes problem) has been the subject of numerous studies, and is known to be consistent up to a time step dependent splitting error which vanishes in the continuous time limit (see, e.g.,
\cite{Marchuk1990,Guermond2006,Liu2010} 
among many others). However, specific to a particle-mesh framework, we also have to require consistency of the projection operators $\mathcal{P}_E$ (particle - mesh) and $\mathcal{P}_L$ (mesh - particle). 
This implies that subsequent application of the projection operators $\mathcal{P}_L$ and $\mathcal{P}_E$ must recover the initial scalar-valued mesh field $\psi_h$ exactly. Similarly, provided that the particle field $\psi_p$ is initialized by means of the mesh-particle projection $\mathcal{P}_L$, subsequent application of the projection operators $\mathcal{P}_E$ and $\mathcal{P}_L$ must recover the initial particle field $\psi_p$ exactly. Mathematically, these compatibility conditions can be formulated as, respectively
    \begin{subequations} \label{eq:compatibility_projections}
    \begin{align}
    \mathcal{P}_E \circ \mathcal{P}_L(\xDiscreteScalar{\psi}) & = \xDiscreteScalar{\psi}, \\
    \mathcal{P}_L \circ \mathcal{P}_E(\psi_p) & = \psi_p,
    \end{align}
    \end{subequations}
which implies mutual consistency of steps 2 and 4 in the absence of steps 1 and 3.
To verify these conditions for the presented formulations for $\mathcal{P}_E$ and $\mathcal{P}_L$, consider the scalar-valued projections. For a repeated back-and-forth mapping between particles and mesh, hence omitting step 1 (i.e. $\mathbf{a}=0$) and step 3 (i.e. $\kappa =0$), the constraint will be inactive in Eq.~\eqref{eq:scalar_optimality_conditions}, i.e. $\xDiscreteScalar{\lambda} = 0$ everywhere, and $\xDiscreteScalar{\bar{\psi}}$ and $\xDiscreteScalar{\psi}$ coincide in a weak sense at cell facets. As a result, the only term remaining in the particle-mesh projection stems from the $\ell^2$-term in the Lagrangian functional $\mathcal{L}$ given by Eq.~\eqref{eq:scalar_lagrangian-functional}. Hence, the functional $\mathcal{L}$ underpinning the particle-mesh projection is in this case similar to the objective function $J$ in Eq.~\eqref{eq:J1}. By recalling that a trivial projection between $\xDiscreteScalar{\psi} \in \xDiscreteScalar{W} $ and $\xDiscreteScalar{\phi} \in \xDiscreteScalar{W}$ can be established, it follows that the mesh-particle projection (Eq.~\eqref{eq:mesh2particle_semidiscrete}) is based on the same objective function, and only differs in terms of the minimizer.
Owing to this symmetry Eq.~\eqref{eq:compatibility_projections} is satisfied, which shows consistency of the projections in the semi-discrete setting.

Concerning (ii), consistency of the HDG method used in the diffusion step was proven in \cite{Labeur2007}. For a consistency proof of the HDG discretization of the Stokes equations, reference is made to \cite{Labeur2012,Rhebergen2016}.

It remains to prove consistency of the PDE-constrained projection (iii). To this end, consider a sufficiently smooth scalar field $\psi$. Substitution into Eq.~\eqref{eq:semi_discr_state_eq} gives, after integration by parts,
    \begin{multline}
    	 \sum_{K} \areaIntegral{\Omega}{\brac{\frac{\partial \psi}{\partial t} + \nabla \cdot \brac{\mathbf{a} \psi}} \delta \lambda_h}
    	+
    	\sum_{K}^{}\lineIntegral{\partial K\setminus \varGamma_N}{ \mathbf{a} \cdot \mathbf{n}  \, \brac{\bar{\psi}_h - \psi_h}\delta \lambda_h } \\
    	- 
    	\lineIntegral{\varGamma_N}{\gamma \mathbf{a}\cdot\mathbf{n} \, \psi \; \delta \lambda_h} 
            +   
    	\lineIntegral{\varGamma_N}{\gamma h_a \; \delta \lambda_h} 
    	=
    	0 \hspace{10pt} \quad \forall \; \delta \lambda_h \in T_h,
    \end{multline}
which demonstrates consistency with the strong form of the advection problem  Eqs.~(\ref{eq:scalar_hyperbolic_subproblem_transport}-\ref{eq:scalar_hyperbolic_subproblem_flux}) and the Neumann boundary condition Eq.~\eqref{eq:scalar_hyperbolic_subproblem_neumann}, with the enforcement of $\bar{\psi}=\psi$ on interior cell facets and the Dirichlet boundary $\Gamma_D$. 
\subsubsection{Conservation}
We will prove that the scheme is mass and momentum conservative from the perspective of the Eulerian mesh. Global and local conservation of the HDG method for the Eulerian diffusion step was demonstrated in \cite{Labeur2007}. Global and local conservation of the HDG method for the unsteady Stokes step was demonstrated in \cite{Labeur2012,Rhebergen2016,Rhebergen2017}. It therefore remains to prove conservation of the PDE-constrained particle-mesh projection. 

Setting $\delta \lambda_h=1$ in Eq.~\eqref{eq:semi_discr_state_eq}, and rearranging, leads to
\begin{equation} \label{eq:global_conservation}
\areaIntegral{\Omega}{\frac{\partial \psi_h }{\partial t}}    =
- \sum_{K}^{}\lineIntegral{\partial K\setminus \varGamma_N}{ \mathbf{a} \cdot \mathbf{n}  \, \bar{\psi}_h } 
- \lineIntegral{\varGamma_N}{\brac{1-\gamma} \mathbf{a} \cdot \mathbf{n} \, \psi_h}  
- \lineIntegral{\varGamma_N}{\gamma h_a }.  
\end{equation}
For a point-wise divergence free vector field $\mathbf{a}$ the boundary integral on the union of interior cell facets vanishes, due to $\bar{\psi}_h$ being single-valued on facets $F \in \mathcal{F}$. The right-hand side therefore equals the total ingoing advective flux at the exterior boundary $\varGamma$, thereby proving global mass conservation. 

For local mass conservation, setting $\delta \lambda_h = 1$ on cell $K$ and $\delta \lambda_h=0$
on $\Omega \setminus K$ gives, after rearrangement,
\begin{equation} \label{eq:local_conservation}
	 \areaIntegral{K}{\frac{\partial \psi_h }{\partial t}}    =
	- \lineIntegral{\partial K\setminus \varGamma_N}{ \mathbf{a} \cdot \mathbf{n}  \, \bar{\psi}_h } 
         - \lineIntegral{\varGamma_N}{\brac{1-\gamma} \mathbf{a} \cdot \mathbf{n} \, \psi_h}  
	- \lineIntegral{\varGamma_N}{\gamma h_a }  
\end{equation}
The right-hand side of Eq.~\eqref{eq:local_conservation} constitutes the ingoing advective flux on the cell facet $\partial K$ which proves local conservation in terms of the numerical flux on $\mathcal{F}$.
For the Navier-Stokes problem, noting that the advective field $\xDiscreteVector{\bar{U}} \in \boldsymbol{\bar{W}}_{h,g}$ is single-valued on $\mathcal{F}$, the derivation of the global and local conservation statements for the vector-valued quantity $\xDiscreteVector{v} \in \xDiscreteVector{W}$ in Eq.~\eqref{eq:optimality_NS} proceeds along similar lines. 

In this work particles carry point evaluations of the underlying field with no notion of mass or volume. As such, a notion of conservation is not well-defined at the particle level. Therefore, the current approach can only pursue a notion of conservation at the mesh level.
\subsubsection{Function spaces} \label{sec:function_space_general_remarks}
To conclude this section, two remarks are made to further elucidate the choice for the particular function space definitions in Section~\ref{sec:problem-definitions}. First, the all-important ingredient for performing the PDE-constrained interaction is to require the optimal solution for the state variable $\psi_h$ (for the scalar-valued case) to be in the discontinuous polynomial space $W_h$, with the scalar-valued facet variable $\bar{\psi}_h \in \bar{W}_{h,g}$ acting as the control variable. With both ingredients being naturally provided by the hybridized Discontinuous Galerkin (HDG) approach, we confidently state that the choice for the HDG method is not just another option, but instead is indispensable in providing the necessary optimality control on the solution. This would be impossible if one were to use, e.g., a continuous or discontinuous Galerkin framework. 

The second remark concerns the choice of function spaces for the semi-discrete form of the unsteady Stokes problem, i.e. $\left(\xDiscreteVector{W}, \boldsymbol{\bar{W}}_{h,g}, Q_h, \bar{Q}_h\right)$. For this particular choice, the resulting velocity fields $\xDiscreteVector{u} \in \xDiscreteVector{W}$ are $H(\text{div})$-conforming and pointwise divergence-free \cite{Rhebergen2017}. In the scope of particle-mesh methods, this has a distinct advantage in that particles can be advected through velocity fields being pointwise-divergence free within a cell and having a continuous normal velocity component across cell facets. In \cite{Maljaars2017}, it was shown that satisfying these two criteria is paramount in maintaining a uniform particle distribution over time. We will further exploit this feature in Section~\ref{sec:particle_advection_ns_discrete}.
\section{Fully-discrete formulations} \label{sec:fully_discrete_formulations}
The final step preceding a computer implementation involves discretization of the respective semi-discrete model components in time. The particular time stepping strategy presented in this section largely follows the approach used in \cite{Maljaars2017}, albeit special care is required to render the particle-mesh projection compatible with the mesh-particle projection. This section concludes by discussing some algorithmic aspects. In particular, we evaluate the choice for the Lagrange multiplier space and present an efficient solution procedure for the PDE-constrained particle-mesh interaction by using static condensation.
\subsection{Advection-diffusion equation}
\subsubsection{Particle advection}
From Eq.~\eqref{eq5:particle_vel_update} it follows that particle quantities other than the position remain constant throughout this stage. A fully-discrete implementation of the Lagrangian advection stage is therefore obtained by integrating Eq.~\eqref{eq5:particle_advection} in time to advance the particle position from $\mathbf{x}_p^n \rightarrow \mathbf{x}_p^{n+1}$. For this purpose, we use a three-stage third-order accurate Runge-Kutta scheme for the linear advection-diffusion problem in which the advective fields are known \textit{a priori}.

At inflow boundaries, particles have to enter the domain. To this end, particles are seeded in the cell contiguous to the inflow boundaries, in such a way to  keep the number of particles constant in these cells. In order to remain consistent with the boundary conditions in Eq.~\eqref{eq:scalar_hyperbolic_subproblem}, properties of the inserted particles are interpolated from the corresponding values imposed by the boundary conditions. 
\subsubsection{PDE-constrained particle-mesh projection} \label{sec:particle-mesh_discrete}
A fully-discrete PDE-constrained particle-mesh projection is formulated with the objective to optimize the scalar field $\psi_h$ at time level $n+1$ given the particle field $\psi_p$ at time $t^{n+1}$. Employing the $\theta$-method, where $1/2 \leq \theta \leq 1$, the constraint in the Lagrangian functional $\mathcal{L}$ is evaluated at time $t^{n+\theta}:= (1-\theta) t^n + \theta t^{n+1}$ using linear interpolation between discrete time levels. The scalar field $\psi_h$ at time level $n+\theta$ is then approximated by
   \begin{equation} \label{eq:interpolated_field}
   \psi_h\left(t^{n+\theta}\right) \approx (1-\theta) \psi_h^{*,n} + \theta \psi_h^{n+1}
   \end{equation}
in which $\psi_h^{*,n} \in W_h$ is an initial field given by
	\begin{equation} \label{eq:consistency_term}
	\xtDiscreteScalar{\psi}{*,n} =  \xtDiscreteScalar{\psi}{n} + \Delta t_n \left((1-\theta_L) \xtDiscreteScalar{\dot{\phi}}{n - 1} + \theta_L  \xtDiscreteScalar{\dot{\phi}}{n} \right),  
	\end{equation}
where $\Delta t_n = t^{n+1} - t^n$ is the time step size, $\theta_L$ is an additional time stepping parameter ($1/2 \leq \theta_L \leq 1$, but possibly different from $\theta$), and the increments $\dot{\phi}_h^m$ (with $m=n-1,n$) are defined by
\begin{equation} \label{eq:mesh_increment}
    \dot{\phi}_h^m = \frac{\phi_h^m - \phi_h^{*,m-1}}{\Delta t_{m-1}},
\end{equation}
with $\phi_h^m$ and $\phi_h^{*,m}$ being fields in $W_h$ related to the Eulerian diffusion step, to be formulated shortly. Using these increments in Eq.~\eqref{eq:consistency_term}, insures compatibility of the fully-discrete projection operator $\mathcal{P}_L$ with $\mathcal{P}_E$ which is required for consistency, as discussed in Section~\ref{sec:consistency1} and to be elaborated further in Section~\ref{sec:time_accuracy_considerations}. The time derivative of the scalar field $\psi_h$ at time level $n+\theta$ is now given by
\begin{equation} \label{eq:discr_time_derivative}
\left. \partialDt{\xDiscreteScalar{\psi}} \right|_{t^{n+\theta}}\approx \frac{\xtDiscreteScalar{\psi}{n+1} - \xtDiscreteScalar{\psi}{*,n}}{\Delta t_n },
\end{equation}
which follows from the linear interpolation used in Eq.~\eqref{eq:interpolated_field}. \\
Next, variations of the dependent fields are taken with respect to the degrees of freedom at time level $n+1$, which involves the replacement of variations $\left(\delta \xDiscreteScalar{\phi}, \delta \xDiscreteScalar{\lambda}, \delta \xDiscreteScalar{\bar{\phi}} \right) \in \left(\xDiscreteScalar{W}, \xDiscreteScalar{T}, \xDiscreteScalar{\bar{W}}\right)$ in the optimality system Eq.~\eqref{eq:scalar_optimality_conditions} with test functions $\left(\xDiscreteScalar{w}, \xDiscreteScalar{\tau}, \xDiscreteScalar{\bar{w}} \right) \in \left(\xDiscreteScalar{W}, \xDiscreteScalar{T}, \xDiscreteScalar{\bar{W}}\right)$. Using the expression for the time derivative of $\psi_h$ given in Eq.~\eqref{eq:discr_time_derivative}, the time derivative appearing in the co-state equation~\eqref{eq:objective_semi_discrete} is approximated as follows
	\begin{equation} 
		\left. \frac{\partial \delta \psi_h}{\partial t} \right|_{t^{n+\theta}} 
		\approx 
		\frac{\delta \psi_h^{n+1} - \delta \psi_h^{*,n}}{\Delta t_n} = \frac{w_h}{\Delta t_n},
	\end{equation} 
since variations $\delta \xtDiscreteScalar{\psi}{*,n} \in W_h$ vanish.

Given these approximations, the fully-discrete co-state equation reads: given the particle field $\psi_p^{n} \in \Psi_t$, the post-advection particle positions $\mathbf{x}_p^{n+1} \in \mathcal{X}_t$, and the intermediate field $\xtDiscreteScalar{\psi}{*,n}\in W_h$, find $\left(\xtDiscreteScalar{\psi}{n+1}, \xtDiscreteScalar{\lambda}{n+1}, \xtDiscreteScalar{\bar{\psi}}{n+1} \right) \in  \left( W_h, T_h, \bar{W}_{h,g} \right)$ such that
\begin{subequations} \label{eq:discrete_optimality-adv-diff}
\begin{multline}\label{eq:costate-discrete-adv-diff}
\sum_{p \in \mathcal{S}_t} \left( \xtDiscreteScalar{\psi}{n+1}(\mathbf{x}_p^{n+1}) - 
\psi_p^{n} \right) w_h(\mathbf{x}_p^{n+1})   
-
\sum\limits_{K}^{} \lineIntegral{\partial K}{\beta \left( \xtDiscreteScalar{\bar{\psi}}{n+1} - \xtDiscreteScalar{\psi}{n+1} \right) \xDiscreteScalar{w}} \\
+ 
\areaIntegral{\Omega}{\frac{\xDiscreteScalar{w}}{\Delta t_n} \lambda^{n+1}_h } 
- 
\theta \sum_{K}\areaIntegral{K}{ \left(\mathbf{a}^{} \xDiscreteScalar{w} \right) \cdot \nabla \xtDiscreteScalar{\lambda}{n+1} }  
+
\theta \lineIntegral{\NeumannBC}{ \left(1 - \gamma \right) \mathbf{a} \cdot \mathbf{n} \xtDiscreteScalar{\lambda}{n+1}  \xDiscreteScalar{w}}
= 0 \hspace{10pt} \forall \xDiscreteScalar{w} \in W_h.
\end{multline}
Correspondingly, the fully-discrete counterpart of the state equation Eq.~\eqref{eq:semi_discr_state_eq} reads:
\begin{multline}\label{eq:state-discrete-adv-diff}
\areaIntegral{\Omega}{\frac{\xtDiscreteScalar{\psi}{n+1} - \xtDiscreteScalar{\psi}{*,n} }{\Delta t_n} \xDiscreteScalar{\tau}  }    
- 
\theta \sum_{K}\areaIntegral{K}{\left(\mathbf{a} \xtDiscreteScalar{\psi}{n+1}\right)  \cdot \nabla \xDiscreteScalar{\tau} }  
+
\sum_{K}^{}\lineIntegral{\partial K\setminus \NeumannBC}{\mathbf{a} \cdot \mathbf{n} \xtDiscreteScalar{\bar{\psi}}{n+1} \xDiscreteScalar{\tau}} 
+
\theta \lineIntegral{\NeumannBC}{\left(1-\gamma \right) \mathbf{a} \cdot \mathbf{n} \xtDiscreteScalar{\psi}{n+1} \xDiscreteScalar{\tau}} \\
+
\lineIntegral{\NeumannBC}{\gamma h^{n+\theta}_a \xDiscreteScalar{\tau}}
=
(1-\theta) \sum_{K}\areaIntegral{K}{ \left(\mathbf{a} \xtDiscreteScalar{\psi}{*,n}\right)  \cdot \GradientScalar{\xDiscreteScalar{\tau}} }
-
(1-\theta) \lineIntegral{\NeumannBC}{\left(1-\gamma \right) \mathbf{a} \cdot \mathbf{n} \xtDiscreteScalar{\psi}{*,n} \xDiscreteScalar{\tau}} 
 \hspace{10pt} \forall \xDiscreteScalar{\tau} \in T_h.
\end{multline}
Finally, the fully-discrete optimality condition becomes
\begin{align}
\sum_{K}^{}\lineIntegral{\partial K\setminus\varGamma_N}{ \mathbf{a} \cdot \mathbf{n} \xtDiscreteScalar{\lambda}{n+1} \bar{w}_h }  
+
\sum\limits_{K}^{} \lineIntegral{\partial K}{\beta \left( \xtDiscreteScalar{\bar{\psi}}{n+1} - \xtDiscreteScalar{\psi}{n+1} \right) \bar{w}_h }
&= 0 
&& \forall \xDiscreteScalar{\bar{w}} \in \bar{W}_{h,0} \label{eq:optimality-discrete-adv-diff}.
\end{align}
\end{subequations}
In these equations, The Lagrange multiplier $\lambda_h$ and the control variable $\Bar{\psi}_h$ are conveniently taken at time level $n+1$, which is allowed since these variables are fully-implicit, not requiring differentiation in time. 

The reconstructed field $\xtDiscreteScalar{\psi}{n+1}$
which is obtained after solving Eq.~\eqref{eq:discrete_optimality-adv-diff} for $(\xtDiscreteScalar{\psi}{n+1}, \xtDiscreteScalar{\lambda}{n+1}, \xtDiscreteScalar{\bar{\psi}}{n+1} )$ 
will serve as an input to the solver for the diffusion equation. 
\subsubsection{Diffusion equation} \label{sec:diffusion_fully_discrete}
A fully-discrete formulation for the HDG discretization of the diffusion equation (Eq.~\eqref{eq:diffusion_semi_discrete}) is formulated using the backward Euler method to discretize in time. The resulting discrete problem for the respective components then reads: given the initial field $\xtDiscreteScalar{\phi}{*,n} = \xtDiscreteScalar{\psi}{n+1}  \in W_h$, the source term $f^{n+1}$, the diffusive Neumann boundary condition $h_d^{n+1}$ and the diffusivity $\kappa$, find $\xtDiscreteScalar{\phi}{n+1} \in W_h$ and $\xtDiscreteScalar{\bar{\phi}}{n+1} \in \bar{W}_{h,g}$ such that 
	\begin{subequations} \label{eq:diffusion_discrete}
	\begin{equation}
	\begin{alignedat}{2}
	\areaIntegral{\Omega}{\frac{\xtDiscreteScalar{\phi}{n+1} - \xtDiscreteScalar{\phi}{*,n}}{\Delta t_n} \xDiscreteScalar{w}} 
	+
	\sum\limits_{K}^{}\areaIntegral{K}{\kappa \nabla \xtDiscreteScalar{\phi}{n+1} \cdot \nabla \xDiscreteScalar{w}} 
    +
	\sum\limits_{K}^{}\lineIntegral{\partial K}{ \boldsymbol{\hat{\sigma}}_{d,h}^{n+1} \cdot \mathbf{n} \xDiscreteScalar{w}} 
	\\
	+
	\sum\limits_{K}^{}\lineIntegral{\partial K}{\kappa \left(\xtDiscreteScalar{\bar{\phi}}{n+1} - \xtDiscreteScalar{\phi}{n+1} \right)\mathbf{n} \cdot \nabla \xDiscreteScalar{w}}
	= 
	\areaIntegral{\Omega}{f^{n+1} \xDiscreteScalar{w}} && \qquad \forall \hspace{3pt} \xDiscreteScalar{w} \in W_h,
	\end{alignedat}
	\end{equation}
	and 
	\begin{align}
	\sum\limits_{K}^{}\lineIntegral{\partial K}{ \boldsymbol{\hat{\sigma}}_{d,h}^{n+1} \cdot \mathbf{n} \xDiscreteScalar{\bar{w}}} 
	= 
	\lineIntegral{\NeumannBC}{h_d^{n+1} \xDiscreteScalar{\bar{w}}} && \forall \hspace{3pt} \xDiscreteScalar{\bar{w}}\in \bar{W}_{h,0},
	\end{align}
	in which the fully-discrete diffusive flux $\boldsymbol{\hat{\sigma}}_{d,h}^{n+1}$ is given by
	\begin{equation}
	    \boldsymbol{\hat{\sigma}}_{d,h}^{n+1} = -\kappa \nabla \xtDiscreteScalar{\phi}{n+1}- \frac{\alpha}{h_K} \kappa  \left(\xtDiscreteScalar{\xtDiscreteScalar{\Bar{\phi}}{n+1} - \phi}{n+1}\right).
	\end{equation}
	\end{subequations}
Solving Eq.~\eqref{eq:diffusion_discrete} for $\left(\xtDiscreteScalar{\phi}{n+1}, \xtDiscreteScalar{\bar{\phi}}{n+1}\right)$ provides the ingredients necessary to update the particle field $\psi_p$ in the subsequent mesh-particle projection step. 
\subsubsection{Mesh-particle projection} \label{sec:mesh2particle_discrete}
Finally, a fully-discrete mesh-particle projection is formulated by mapping the increments of the mesh related field ($\phi_h$) to the particles using
\begin{align}\label{eq:discrete_particle_update}
\psi_p^{n+1} = \psi_p^n + \Delta t_n \left( (1- \theta_L) \xtDiscreteScalar{\dot{\phi}}{n}\left(\mathbf{x}_p^n\right) + \theta_L \xtDiscreteScalar{\dot{\phi}}{n+1}\left(\mathbf{x}_p^{n+1}\right) \right) && \forall \hspace{3pt} p \in \mathcal{S}_t,
\end{align}
where $1/2 \leq \theta_L \leq 1$, and $\xtDiscreteScalar{\dot{\phi}}{n} \in W_h$ is defined according to Eq.~\eqref{eq:mesh_increment}. 
We emphasize the consistency of the formulation for the mesh-particle update (Eq.~\eqref{eq:discrete_particle_update}) with the definition of $\xtDiscreteScalar{\psi}{*,n} \in W_h$ in Eq.~\eqref{eq:consistency_term}, which is required to respect the compatibility condition, Eq.~\eqref{eq:compatibility_projections}, in the fully-discrete setting. 
Furthermore, it readily follows that in the advective limit (i.e. for $\kappa = 0$) it holds that $\psi_p^{n+1} = \psi_p^n$, since $\dot{\phi}^{n}_{h} = \dot{\phi}^{n+1}_{h} = 0$. 
\subsection{A fully-discrete formulation for the incompressible Navier-Stokes equations} \label{sec:navier_stokes_fully_discrete}
We now briefly present the fully-discrete formulations for the four model components constituting the particle-mesh operator splitting method for the incompressible Navier-Stokes equations. The formulations are quite similar to those for the advection-diffusion problem but, in addition, involve the incompressibility constraint and the discretization of the advective fields $\xDiscreteVector{U} \in \xDiscreteVector{W}$ and $\xDiscreteVector{\bar{U}} \in \boldsymbol{\bar{W}}_{h,g}$.
\subsubsection{Particle advection} \label{sec:particle_advection_ns_discrete}
Once the advective velocity field $\xDiscreteVector{U} \in \xDiscreteVector{W}$ has been chosen, the particle advection for the incompressible Navier-Stokes equations proceeds along similar lines as for the linear advection-diffusion equation. To this end, a judicious choice for the advective field $\xDiscreteVector{U}(\mathbf{x},t)$ in Eq.~\eqref{eq5:ns_particle_advection} for $t \in I_n$, where $I_n =(t^n, t^{n+1}]$, is to use the mesh related velocity field $\mathbf{u}_h \in \mathbf{W}_h$ frozen at the old time level $t^n$, i.e.
\begin{equation} \label{eq:advective_field}
\xDiscreteVector{U}\left(\mathbf{x}_p(t),t\right) = \mathbf{u}_h(\mathbf{x}_p(t),t^n).
\end{equation}

The $H(\text{div})$-conformity of the velocity field $\mathbf{u}_h \in \mathbf{W}_h$ following the solution of the Stokes problem Eq.~\eqref{eq:unsteady-stokes_semi-discrete} guarantees that the advective field $\xDiscreteVector{U}$ is pointwise divergence-free (see Section~\ref{sec:function_space_general_remarks}), while it is also explicitly known in each time interval $I_n \in \mathcal{I}$. Thus, Eq.~\eqref{eq5:ns_particle_advection} can be integrated in time using explicit integration schemes. 
To this end, a second-order two-step Adams-Bashforth scheme is used. This choice is motivated by considering the time accuracy of the discretization of the Stokes equation on the mesh in conjunction with the projection steps, see Sections~\ref{sec:ns_fully_discrete_pde}~-~\ref{sec:ns_disrete_mesh_particle}. 
The combination of these steps can be rendered second-order accurate in time, as will be argued in Section~\ref{sec:time_accuracy_considerations}.
\subsubsection{PDE-constrained particle-mesh projection} \label{sec:ns_fully_discrete_pde}
The fully-discrete vector-valued PDE-constrained particle-mesh projection  is formulated similarly to its scalar-valued counterpart in that variations of the dependent fields are considered with respect to the degrees of freedom at time level $n+1$. This involves the replacement of variations $\left( \delta \xDiscreteVector{v}, \delta \xDiscreteVector{\lambda}, \delta \xDiscreteVector{\bar{v}} \right)$ in Eq.~\eqref{eq:optimality_NS} with test functions $\left( \xDiscreteVector{w}, \xDiscreteVector{\tau}, \xDiscreteVector{\bar{w}} \right) \in \left( \xDiscreteVector{W}, \xDiscreteVector{T}, \xDiscreteVector{\bar{W}}\right)$.
Once more, the $\theta$-method (with $1/2\leq\theta\leq1$) is used to integrate in time with the optimality system Eq.~\eqref{eq:optimality_NS} being evaluated at time $t^{n+\theta}$ using linear interpolation between time levels $n$ and $n+1$. The time derivative appearing in the state equation Eq.~\eqref{eq:state-semi-discrete_NS} is approximated as follows 
\begin{equation}
\left. \partialDt{\xDiscreteVector{v}} \right|_{t^{n+\theta}} \approx \frac{\xtDiscreteVector{v}{n+1} -  \xtDiscreteVector{v}{*,n}}{\Delta t_n }, 
\end{equation}
in which, analogously to Eq.~\eqref{eq:consistency_term}, $\xtDiscreteVector{v}{*,n}$  is defined as 
\begin{equation} \label{eq:old_tlevel_consisten_NS}
\xtDiscreteVector{v}{*,n} := \xtDiscreteVector{v}{n} + \Delta t_n \left( (1-\theta_L)  \Dot{\mathbf{u}}^{n-1} + \theta_L \Dot{\mathbf{u}}^{n} \right),
\end{equation}
where the acceleration field $\Dot{\mathbf{u}}^{n} \in \xDiscreteVector{W}$ (i.e., the vector-valued counterpart of Eq.~\eqref{eq:mesh_increment}) results from the fully-discrete Stokes step (to be discussed in Section~\ref{sec:stokes_discrete}). Consistent with Eq.~\eqref{eq:advective_field}, the advective fields $\xDiscreteVector{U}$ and $\xDiscreteVector{\bar{U}}$ in Eqs.\eqref{eq:optimality_NS} are given explicitly by the corresponding velocity fields $\xDiscreteVector{u}$, respectively, $\xDiscreteVector{\bar{u}}$ from the previous 
Stokes step, i.e. for $t \in I_n$,
\begin{align}
\xDiscreteVector{U}(\mathbf{x}, t)  	:=  \xtDiscreteVector{u}{n}(\mathbf{x}), && \text{and} && \xDiscreteVector{\bar{U}}(\mathbf{x}, t)  	:=  \xtDiscreteVector{\bar{u}}{n}(\mathbf{x}),
\end{align}
which \textit{de facto} linearizes the PDE-constrained particle-mesh projection problem.

Given these definitions, we now present a fully-discrete implementation of the optimality system Eq.~\eqref{eq:optimality_NS}. 
The time-discrete counterpart of the co-state equation (Eq.~\eqref{eq:costate-semi-discrete_NS}) reads as follows: 
given the vector-valued particle field $\ContVector{v}_p^{n} \in \mathcal{V}_t$, the particle positions $\mathbf{x}_p^{n+1} \in \mathcal{X}_t$, and the field $\xtDiscreteVector{v}{*,n} \in \xDiscreteVector{W}$, find $\left( \xtDiscreteVector{v}{n+1}, \xtDiscreteVector{\lambda}{n+1}, \xtDiscreteVector{\bar{v}}{n+1}\right) \in \left( \xDiscreteVector{W}, \xDiscreteVector{T}, \boldsymbol{\bar{W}}_{h,g} \right)$ such that
\begin{subequations} \label{eq:optimality_NS_discrete}
\begin{multline}\label{eq:costate-discrete_NS}
\sum_{p \in \mathcal{S}_t} \left( \xtDiscreteVector{v}{n+1}(\mathbf{x}_p^{n+1}) - \mathbf{v}_p^{n} \right) \cdot \xDiscreteVector{w}(\mathbf{x}_p^{n+1})   
- 
\sum\limits_{K}^{} \lineIntegral{\partial K}{\beta \left( \xtDiscreteVector{\bar{v}}{n+1} - \xtDiscreteVector{v}{n+1} \right) \cdot \xDiscreteVector{w}}
+
\areaIntegral{\Omega}{\frac{\xDiscreteVector{w}}{\Delta t_n} \cdot \xtDiscreteVector{\lambda}{n+1} }   \\
- 
\theta \sum_{K}\areaIntegral{K}{ \left(\xDiscreteVector{w} \otimes \xtDiscreteVector{u}{n} \right) \colon \GradientVector{\lambda}^{n+1} }  
+
\theta \lineIntegral{\NeumannBC}{(1-\gamma) \left(\xDiscreteVector{w} \otimes \xtDiscreteVector{u}{n}\right) \mathbf{n} \cdot \xtDiscreteVector{\lambda}{n+1}}
= 
0 \hspace{10pt} \forall \hspace{3pt}  \xDiscreteVector{w} \in \xDiscreteVector{W}.
\end{multline}
Simultaneously the fully-discrete state equation (Eq.~\eqref{eq:state-semi-discrete_NS}) is to be satisfied, given the advective Neumann boundary condition $\mathbf{h}_a$,
\begin{multline} \label{eq:state-discrete_NS}
\areaIntegral{\Omega}{\frac{\xtDiscreteVector{v}{n+1} - \xtDiscreteVector{v}{*,n}  }{\Delta t_n} \cdot \xDiscreteVector{\tau}  }    
- 
\theta \sum_{K}\areaIntegral{K}{\left( \xtDiscreteVector{v}{n+1} \otimes  \mathbf{u}^{n}_h\right)  : \nabla \xDiscreteVector{\tau} }  
+ 
\sum_{K}^{}\lineIntegral{\partial K\setminus \NeumannBC}{\left(\xtDiscreteVector{\bar{v}}{n+1} \otimes \xtDiscreteVector{\bar{u}}{n} \right) \mathbf{n} \cdot \xDiscreteVector{\tau}}  \\
+
\theta \lineIntegral{\NeumannBC}{(1-\gamma) \left(\xtDiscreteVector{v}{n+1} \otimes \xtDiscreteVector{u}{n}\right) \mathbf{n} \cdot \xDiscreteVector{\tau}}
+ 
\lineIntegral{\NeumannBC}{ \gamma \mathbf{h}_a^{n+\theta} \cdot \xDiscreteVector{\tau}} 
\\  
=
(1-\theta) \sum_{K}\areaIntegral{K}{ \left( \xtDiscreteVector{v}{*,n} \otimes  \xtDiscreteVector{u}{n} \right)  : \nabla \xDiscreteVector{\tau} }
-
(1-\theta) \lineIntegral{\NeumannBC}{(1-\gamma) \left(\xtDiscreteVector{v}{*,n} \otimes \xtDiscreteVector{u}{n}\right) \mathbf{n} \cdot \xDiscreteVector{\tau}}
 \hspace{10pt} \forall \hspace{3pt}  \xDiscreteVector{\tau} \in \xDiscreteVector{T},
\end{multline}
together with the fully-discrete optimality condition (Eq.~\eqref{eq:optimality-semi-discrete_NS}), which reads,
\begin{equation}\label{eq:optimality-discrete_NS}
\sum_{K}^{}\lineIntegral{\partial K\setminus \NeumannBC}{\left( \xtDiscreteVector{\lambda}{n+1} \otimes \xtDiscreteVector{u}{n} \right) \mathbf{n} \cdot \xDiscreteVector{\bar{w}}}   
+
\sum\limits_{K}^{} \lineIntegral{\partial K}{\beta \left( \xtDiscreteVector{\bar{v}}{n+1} - \xtDiscreteVector{v}{n+1} \right) \cdot \xDiscreteVector{\bar{w}} }
= 0 
\hspace{10pt} \forall \hspace{3pt} \mathbf{\bar{w}}_h \in \boldsymbol{\bar{W}}_{h,0}. 
\end{equation}
\end{subequations}
Solving Eqs.~\eqref{eq:optimality_NS_discrete} for $\left( \xtDiscreteVector{v}{n+1}, \xtDiscreteVector{\lambda}{n+1}, \xtDiscreteVector{\bar{v}}{n+1}\right)$ gives the reconstructed 
field $\xtDiscreteVector{v}{n+1}$, which will serve as an input to the Stokes solver (to be discussed next).
We emphasize that the fully-discrete state equation~\eqref{eq:state-discrete_NS} does not involve upwinding of the (interface) fluxes.
\subsubsection{Stokes problem} \label{sec:stokes_discrete}
The unsteady Stokes equations are discretized in time by using the backward Euler method, similarly to the method used in \cite{Maljaars2017}. The fully-discrete problem then reads as follows:
given the intermediate field $\xtDiscreteScalar{\mathbf{u}}{*,n} = \xtDiscreteVector{v}{n+1} \in \xDiscreteVector{W}$ , the diffusive Neumann boundary condition $\mathbf{h}_d^{n+1} \in \left[L^2\left(\NeumannBC\right)\right]^d$, the forcing term $\mathbf{f}^{n+1}\in \left[L^2\left(\Omega\right)\right]^d$, and the kinematic viscosity $\nu$, find $\left( \mathbf{u}_h^{n+1},\bar{\mathbf{u}}_h^{n+1}, p_h^{n+1}, \bar{p}_h^{n+1} \right) \in \left( \xDiscreteVector{W}, \boldsymbol{\bar{W}}_{h,g}, Q_h, \bar{Q}_h \right) $ such that (local and global momentum balances),
\begin{subequations} \label{eq:stokes_discrete}
\begin{align}
	\begin{split}
		 \areaIntegral{\Omega}{\frac{\xtDiscreteVector{u}{n+1}-\xtDiscreteVector{u}{*,n}}{\Delta t_n}\cdot\xDiscreteVector{w}}
		 -
		 & 
		 \sum\limits_{K}^{}\areaIntegral{K}{\boldsymbol{\sigma}_{d,h}^{n+1} \colon \nabla \xDiscreteVector{w}} 
		 + 
		 \sum\limits_{K}^{}\lineIntegral{\partial K}{\boldsymbol{\hat{\sigma}}_{d,h}^{n+1}\mathbf{n}\cdot \xDiscreteVector{w}} \\ 
		 &
		 +
		 \sum\limits_{K}^{}\lineIntegral{\partial K}{2\nu \left(\bar{\mathbf{u}}_h^{n+1}
		 -
		 \xtDiscreteVector{u}{n+1}\right)\cdot \nabla^s\xDiscreteVector{w}\mathbf{n}} 
		 =  
		 \areaIntegral{\Omega}{\ContVector{f}^{n+1}\cdot \xDiscreteVector{w}} 
		 \hspace{10pt}
		 \forall \hspace{5pt} \xDiscreteVector{w} \in W_h,
	\end{split} \\
	&\sum\limits_{K}^{}\lineIntegral{\partial K}{\boldsymbol{\hat{\sigma}}_{d,h}^{n+1}\mathbf{n} \cdot \xDiscreteVector{\bar{w}} } 
	= 
	\lineIntegral{\varGamma_N}{\ContVector{h}_d^{n+1}\cdot \xDiscreteVector{\bar{w}} }
	\hspace{10pt} \forall \hspace{5pt}
	\xDiscreteVector{\bar{w}} \in \bar{W}_{h,0},
\end{align}
and (local and global mass balances),
\begin{align}
\sum\limits_{K}^{}\areaIntegral{K}{\xtDiscreteVector{u}{n+1}\cdot \nabla q_h} 
- 
\sum\limits_{K}^{}\lineIntegral{\partial K}{\xtDiscreteVector{u}{n+1}\cdot \mathbf{n} \hspace{1pt} q_h} 
&= 
0 \hspace{25pt}
 \forall \hspace{5pt} \xDiscreteScalar{q}\in Q_h,\\
\sum\limits_{K}^{}\lineIntegral{\partial K}{\xtDiscreteVector{u}{n+1}\cdot\mathbf{n}\bar{q}_h} 
- 
\lineIntegral{\partial \Omega}{\xtDiscreteVector{\bar{u}}{n+1}\cdot\mathbf{n}\bar{q}_h} &=0 \hspace{25pt} \forall \hspace{5pt} \xDiscreteScalar{\bar{q}} \in \bar{Q}_h,
\end{align}
\end{subequations}
are satisfied.
In these equations, the diffusive flux $\boldsymbol{\sigma}_{d,h}^{n+1}$ is given by
\begin{equation}
\boldsymbol{\sigma}_{d,h}^{n+1} = \xtDiscreteScalar{p}{n+1} \mathbf{I} - 2\nu\nabla^s \xtDiscreteVector{u}{n+1}, 
\end{equation}
and the diffusive flux at cell facets, $\boldsymbol{\hat{\sigma}}_{d,h}^{n+1}$, is defined as
\begin{equation}
\boldsymbol{\hat{\sigma}}_{d,h}^{n+1} = \xtDiscreteVector{\bar{p}}{n+1}\mathbf{I} - 2\nu \nabla^s \xtDiscreteVector{u}{n+1} - 2\nu \frac{\alpha}{h_K} \left(\xtDiscreteVector{\bar{u}}{n+1} - \xtDiscreteVector{u}{n+1}\right)\otimes \mathbf{n}.
\end{equation}

For a detailed analysis of the unsteady Stokes solver, we refer to \cite{Rhebergen2016,Rhebergen2017}. Here we note that the solver is $\inf-\sup$ stable, conserves momentum globally and locally, is energy stable, produces $H(\text{div})$-conforming velocity fields $\xtDiscreteVector{u}{n+1} \in \xDiscreteVector{W}$, and last but not least, is amenable to a static condensation procedure. 
\subsubsection{Mesh-particle projection} \label{sec:ns_disrete_mesh_particle}
For the incompressible Navier-Stokes equations, the mesh-particle projection is the vector-valued counterpart of Eq.~\eqref{eq:discrete_particle_update}, i.e. the momentum field on particles is updated according to
	\begin{align} \label{eq:crank nicolson}
	\ContVector{v}_p^{n+1} 
	= 
	\ContVector{v}_p^{n} 
	+ 
	\Delta t_n \left((1-\theta_L) \hspace{2pt} \Dot{\mathbf{u}}_h^{n}(\mathbf{x}_p^n) + \theta_L \Dot{\mathbf{u}}_h^{n+1}(\mathbf{x}_p^{n+1}) \right) 
	\hspace{25pt} \forall \hspace{3pt} p \in \mathcal{S}_t,
	\end{align}
with $\Dot{\mathbf{u}}_h^{n}(\mathbf{x}_p^n) $ and $ \Dot{\mathbf{u}}_h^{n+1}(\mathbf{x}_p^{n+1} )$ denoting the mesh-based accelerations at the respective time levels $n$ and $n+1$ - computed according to Eq.~\eqref{eq:mesh_increment} -  which are evaluated at the individual particle positions.
\subsection{Algorithmic aspects} \label{sec:alg_aspects}
\subsubsection{Choosing the Lagrange multiplier space} \label{sec:function_space_considerations_revisited}
A significant simplification of the fully-discrete optimality system is obtained when choosing the polynomial basis for the Lagrange multiplier of the lowest possible order, i.e. $l=0$. For this particular choice, the cell integrals in Eqs.~\eqref{eq:discrete_optimality-adv-diff} and~\eqref{eq:optimality_NS_discrete} containing gradients of the Lagrange multiplier, or of its associated test function, vanish.

This choice will not affect the conservation proofs given by Eqs.~\eqref{eq:global_conservation}~and~\eqref{eq:local_conservation}. Indeed, for the choice $l=0$ it readily follows from Eq.~\eqref{eq:state-discrete-adv-diff} and Eq.~\eqref{eq:state-discrete_NS} that the cell-integrated temporal increments of the conserved variables remain in balance with the facet fluxes. Except for the outflow Neumann boundaries, these facet fluxes are implicitly controlled by the interface variable $\bar{\psi}_h^{n+1} \in \bar{W}_{h,g}$ ($\xtDiscreteVector{\bar{v}}{n+1} \in \boldsymbol{\bar{W}}_{h,g}$), not requiring differentiation in time. Hence, the major advantage of choosing $l=0$ is that it avoids time stepping dependencies in the PDE-constrained particle-mesh projection, thus rendering the particle-mesh projection independent of $\theta$. 

Exploiting this evident advantage, we will restrict the discussion in the sequel to the choice $l=0$. In fact, it is only for this particular choice of the Lagrange multiplier space that second-order accuracy in time can be expected for the given time integration scheme. 
Furthermore, independent of the choice for $l$, we expect optimal spatial convergence rates of the scheme of order $k+1$, provided a sufficiently accurate particle advection scheme is used. 
\subsubsection{Accuracy of the time stepping scheme} \label{sec:time_accuracy_considerations}
For an overview, let us briefly summarize the essential ingredients required for obtaining second-order accuracy in time. In \cite{Maljaars2017}~-~with the particle-mesh projection being formulated in terms of local $\ell^2$-projections~-~it was shown that overall second-order accuracy in time can be expected when combining a backward Euler method for the Eulerian step, Eq.~\eqref{eq:diffusion_discrete} or Eq.~\eqref{eq:stokes_discrete}, with a second-order mesh-particle update as in Eq.~\eqref{eq:discrete_particle_update} or Eq.~\eqref{eq:crank nicolson}, using $\theta_L = 1/2$.  Note that the latter requires saving an earlier particle field. 
Backward Euler suffices on the Eulerian step since $\phi_h^{n+1}$ and $u_h^{n+1}$ are advanced over one time step from the second-order fields $\phi_h^{*,n}$ and $u_h^{*,n}$, respectively. Thus, the local, second-order accuracy of backward Euler pertains.

To extend this approach to the PDE-constrained particle-mesh projection pursued here, we have to use $l=0$ as argued above. In addition, we have to ensure compatibility between the mesh-particle update, Eq.~\eqref{eq:discrete_particle_update} or Eq.~\eqref{eq:crank nicolson}, and the PDE-constrained particle-mesh interaction in the fully-discrete setting, Eq.~\eqref{eq:discrete_optimality-adv-diff} or Eq.~\eqref{eq:optimality_NS_discrete}.
As further detailed in \ref{sec:appendix_consistency}, the latter requirement leads to definitions for $\xtDiscreteScalar{\psi}{*,n}$ or $\xtDiscreteVector{v}{*,n}$ via Eq.~\eqref{eq:consistency_term} or Eq.~\eqref{eq:old_tlevel_consisten_NS}, respectively.
\subsubsection{Numerical implementation} \label{sec:numerical_implementation}
Similarly to the fully-discrete diffusion equation and the fully-discrete Stokes equations, the fully-discrete PDE-constrained particle-mesh projections can be implemented efficiently by eliminating the unknowns local to a cell $K$ (i.e. $\xtDiscreteScalar{\psi}{n+1}$ and $\xtDiscreteScalar{\lambda}{n+1}$ for the scalar-valued case) in favor of the global control variable $\xtDiscreteScalar{\bar{\psi}}{n+1}$ via static condensation. As a result, the size of the global system is reduced significantly as it depends only on the number of degrees of freedom of the control variable. After solving this system for $\xtDiscreteScalar{\bar{\psi}}{n+1}$, the local variables $\xtDiscreteScalar{\psi}{n+1}$ and $\xtDiscreteScalar{\lambda}{n+1}$ are obtained in a back substitution step.
\section{Numerical examples} \label{sec: Numerical Examples}
In this section, the properties of the proposed method are illustrated by means of various numerical examples. 
The first examples illustrate the behavior of the scheme for linear advection-diffusion problems in terms of accuracy and conservation. We consider rotating Gaussian hump, rigid body rotation and skew advection benchmarks to illustrate the method's ability to handle smooth fields, point/line singularities and inflow/outflow conditions, respectively. 
The second part of this section evaluates the performance of the developed approach for the incompressible Navier-Stokes equations by considering plane Poiseuille flow and Taylor-Green flow. 

Throughout, domains $\Omega \in \mathbb{R}^2$ are considered and the time domain of interest is partitioned using constant time step sizes $\Delta t$. Furthermore, the regularization term $\beta$ is set to a fixed value of $10^{-6}$ for all computations. The penalty parameter $\alpha$ in Eq.~\eqref{eq:diffusion_discrete} (diffusion equation) and Eq.~\eqref{eq:stokes_discrete} (Stokes problem) is set to $12k^2$ and $6k^2$, respectively. Unless otherwise specified, we choose $l=0$, thus rendering the scheme independent of $\theta$, see Section~\ref{sec:function_space_considerations_revisited}.

Tools from the finite element framework FEniCS \cite{Logg2012} are used to assemble and solve the equations on the mesh arising from the discretization of the PDE-constrained projections, the diffusion equation and the unsteady Stokes equations. A static condensation procedure is applied in all cases and the resulting global systems are solved using direct Gaussian elimination. 
The computer code for performing the computations (except for the test conducted in Section~\ref{sec:taylor_green_low_reynolds}) is available under an open source license and can be obtained via \href{https://bitbucket.org/jakob_maljaars/leopart}{\texttt{bitbucket.org/jakob\_maljaars/leopart}}.
%
\subsection{Advection-diffusion: Gaussian hump}
The accuracy of the presented method is assessed by considering a (rotating) Gaussian pulse in the diffusive limit, for moderate diffusion, and in the advective limit. 
The domain of interest is the circular disk $\Omega := \{(x,y)\hspace{1pt}\vert\hspace{1pt} x^2+y^2 \leq 0.5  \}$ and the velocity field is either set to $\mathbf{a} = \mathbf{0}$ (diffusive limit), or given by 
    \begin{equation} \label{eq:solid_body_rotation}
    \mathbf{a} = \pi \left(-y, x\right)^\top.
    \end{equation}
The corresponding analytical solution for a rotating Gaussian pulse is given by
	\begin{equation}\label{eq:gaussian_pulse_exact}
	\phi(\mathbf{x},t) 
	= 
	\frac{2 \sigma ^2 }{2 \sigma ^2 + 4 \kappa t} \exp\left(- \frac{\lVert \mathbf{\bar{x}}(\mathbf{x},t) - \mathbf{x}_c \rVert ^ 2  }{ 2 \sigma ^2 + 4 \kappa t} \right), 
	\end{equation}
in which $\mathbf{x}_c = (x_c, y_c)^\top$ is the position vector of the center, $\sigma$ is the initial standard deviation, and $\kappa$ is a constant diffusivity. Furthermore, $ \lVert \cdot \rVert^2$ denotes the square of the Euclidean norm and $\mathbf{\bar{x}}(\mathbf{x},t)$ is a co-ordinate in the rotating frame of reference, given by $\left( x \cos(\pi t) + y \sin(\pi t), -x \sin(\pi t) + y \cos(\pi t)\right)^\top$. 
The initial condition $\phi(\mathbf{x},0)$ is deduced from Eq.~\eqref{eq:gaussian_pulse_exact}, with the initial standard deviation $\sigma$ set to $0.1$ and the Gaussian pulse initially being centered at $(x_c, y_c) = (-0.15, 0)$. The disk-shaped domain $\Omega$ is triangulated using a sequence of mesh refinements, and particles are randomly seeded in $\Omega$ such that each cell contains on average 30 particles, initially. 

Four different cases (listed in Table~\ref{tab:rotating_cone_advectiondiffusion_overview}) are considered for three values of the diffusivity, $\kappa = 0.01$, $\kappa = 0.001$, and $\kappa = 0$ (advective limit). The Dirichlet boundary condition for the diffusion step is deduced from the analytical solution Eq.~\eqref{eq:gaussian_pulse_exact}, and we emphasize that in the advective limit the particle specific mass $\psi_p$ stays constant by virtue of Eq.~\eqref{eq:discrete_particle_update}. The same meshes, initial conditions, and particle distributions are used in all four cases. \\
\begin{table}[H]
\centering
\caption{Gaussian hump: Overview of model settings for advection-diffusion.}
\label{tab:rotating_cone_advectiondiffusion_overview}
\small{
\begin{tabular}{c|c c c c c  }
\hline
\hline 
\rule{0pt}{2.5ex}   
		&	$k$	&	$l$	& $\theta$	& $\theta_L$	& $\mathbf{a}$	\\
\hline 
Case 1	&	2	&	0	& -     &	1/2		&	 $\mathbf{0}$	\\
Case 2	& 	1	& 	0	& -     &	1/2		&	Eq.~\eqref{eq:solid_body_rotation}	\\
Case 3 	& 	2	& 	0	& -     &	1/2		&	Eq.~\eqref{eq:solid_body_rotation}	\\
Case 4 	& 	2	& 	1	& 1/2   &   1/2		&	Eq.~\eqref{eq:solid_body_rotation}	 \\
\hline 
\hline 
\end{tabular}
}
\end{table}
Results obtained after a full revolution are presented in Table~\ref{tab:rotating_cone_advectiondiffusion_errors} for the different test cases. For Cases 1-3, at least second-order convergence is obtained. More precisely, for the largest value of the diffusivity (i.e. $\kappa = 0.01$), the convergence rate tends to second-order, whereas near-optimal convergence is obtained for moderate diffusion, with diffusivity $\kappa = 0.001$, resulting in third-order convergence for Case 1 and Case 3, and second-order convergence for Case 2. In the advective limit, optimal convergence rates are obtained for Case 2 and Case 3. 

Case 3 and Case 4 only differ in the polynomial orders of the Lagrange multiplier space. However, the convergence rates drastically reduce for Case~4 to approximately second-order for the pure advection test (with $\kappa = 0$) and to approximately first order for the mixed advection-diffusion regime. This behavior illustrates the paramount difference between the choice $l=0$ compared to $l\geq1$ in combination with the chosen time stepping scheme, as discussed in Section~\ref{sec:function_space_considerations_revisited}.
\begin{table}[H]
\centering
\caption{Gaussian hump: $L^2$-error for advection-diffusion after one full revolution. Convergence rates based on $h_{K,\text{max}}$.}
\label{tab:rotating_cone_advectiondiffusion_errors}
{\small
\begin{tabular}{c c c c|c c|c c|c c}
\hline
\hline
& & \multicolumn{2}{c|}{Mesh} &  \multicolumn{2}{c|}{$\kappa = 0.01$} &  \multicolumn{2}{c|}{$\kappa = 0.001$} & \multicolumn{2}{c}{$\kappa = 0.0$} \\
& $\Delta t$ & $h_{K,\text{min}}$ & $h_{K,\text{max}}$ & $\lVert \phi - \phi_h \rVert$ & Rate & $\lVert \phi - \phi_h \rVert$ & Rate & $\lVert \phi - \phi_h \rVert$ & Rate  \\ 
\hline
\multirow{ 4}{*}{Case 1} 
&	0.08 &  6.6e-2 	&   1.2e-1 	& 7.3e-5 &  -		& 5.3e-4	& -   & \multicolumn{2}{c}{\multirow{4}{*}{-}} \\
&	0.04 &  3.2e-2 	&   6.2e-2  & 1.5e-5 &  2.3 	& 6.0e-5 	& 3.2 & \multicolumn{2}{c}{} \\
&	0.02 &  1.6e-2 	&   3.1e-2 	& 3.5e-6 &  2.1		& 7.5e-6  	& 3.0 & \multicolumn{2}{c}{} \\
&	0.01 &  7.9e-3 	&   1.6e-2 	& 8.7e-7 &  2.0		& 9.6e-7 	& 3.0 & \multicolumn{2}{c}{} \\
\hline
\multirow{ 4}{*}{Case 2} 
&	0.08  &	6.6e-2 		&   1.2e-1 	& 1.2e-3 	&  - 	& 7.2e-3	&	-	&	1.3e-2	& 	- 	\\	
&	0.04  &  3.2e-2 	&   6.2e-2 	& 3.3e-4 	&  1.9 	& 1.8e-3	&	2.0	&  	3.9e-3	&	1.8	\\
&	0.02  &  1.6e-2 	&   3.1e-2 	& 6.3e-5 	&  2.4 	& 2.7e-4	&	2.7 &	9.6e-4	& 	2.0 \\
&	0.01  &  7.9e-3 	&   1.6e-2 	& 1.2e-5 	&  2.4	& 9.4e-5	&	1.6 & 	2.4e-4	&	2.0	\\
\hline
\multirow{ 4}{*}{Case 3} 
&	0.08  &	6.6e-2 		&   1.2e-1 	& 1.6e-4 	& -		&	8.9e-4	&  - 	&	2.9e-3	&	-	\\
&	0.04  &  3.2e-2 	&   6.2e-2	& 1.6e-5 	& 3.3	&	1.3e-4	& 2.8 	&	2.5e-4	&	3.5	\\
&	0.02  &  1.6e-2 	&   3.1e-2 	& 3.0e-6 	& 2.4 	&	1.9e-5	& 2.8 	&	3.0e-5	&	3.1	\\
&	0.01  &  7.9e-3 	&   1.6e-2 	& 7.4e-7 	& 2.0	&	2.4e-6	& 2.9 	&	4.4e-6	&	2.8	\\
\hline
\multirow{ 4}{*}{Case 4} 
&	0.08  &	 6.6e-2 	&   1.2e-1 	& 3.2e-3	& -		&	1.4e-2	&  - 	&	2.8e-2	&	-	\\
&	0.04  &  3.2e-2 	&   6.2e-2	& 2.0e-3	& 0.7	&	3.4e-3	& 2.0 	&	1.5e-2	&	0.9	\\
&	0.02  &  1.6e-2 	&   3.1e-2 	& 1.1e-3	& 0.9 	&	9.1e-4	& 1.9 	&	3.1e-3	&	2.3	\\
&	0.01  &  7.9e-3 	&   1.6e-2 	& 5.6e-4	& 0.9	&	5.1e-4	& 0.8 	&	8.2e-4	&	1.9	\\
\hline 
\hline  
\end{tabular}
}
\end{table}
\subsection{Advection: rigid body rotation}
In order to qualitatively assess the scheme's ability to preserve point and line singularities, we next consider the advection test proposed in \cite{LeVeque1996}. This test comprises the rigid body rotation of a pointy cone (initially centered at $(x,y)=(-0.3,0)$), a slotted disk (initially centered at $(x,y)=$(0, -0.3)), and a Gaussian hump (initially centered at $(x,y)=$ (0.15, 0.15)) on the the circular disk $\Omega := \{(x,y)\hspace{1pt}\vert\hspace{1pt} x^2+y^2 \leq 0.5  \}$. The velocity field is given by Eq.~\eqref{eq:solid_body_rotation}. By virtue of Eq.~\eqref{eq:discrete_particle_update} it follows that in the advective limit the particle specific masses need not be updated, so that any discontinuities at the particle level persist. 
\begin{table}[H]
\centering
\caption{Rigid body rotation: Overview of mesh, particle and time step settings.}
\label{tab:three_rigid_bodies}
\small{
\begin{tabular}{c c c c c c}
\hline
\hline 
\rule{0pt}{3.ex}   
& $\left|\mathcal{T}\right|$ & $h_{K,\text{min}}$ & $h_{K,\text{max}}$ &$\left|\mathcal{S}_t\right|$ & $\Delta t$ \\
\hline 
Mesh 1  & 16189 & 7.9e-3  &   1.6e-2 &  502480  & 1e-2 \\
Mesh 2 	& 64561 & 4.0e-3  &   7.8e-2 &  2010783 & 5e-3 \\ 
\hline 
\hline 
\end{tabular}
}
\end{table}
A mesh with moderate spatial resolution (Mesh~1, containing 16189 cells), and a fine resolution mesh (Mesh~2, containing 64561 cells) are considered. Approximately 30 particles are assigned per cell initially, resulting in a total number of slightly over $5\times10^5$ and $2\times 10^6$ particles in total for Mesh~1 and Mesh~2, respectively. Time step sizes are chosen such to keep the CFL-number approximately the same on both meshes, see Table~\ref{tab:three_rigid_bodies}.

\begin{figure}[H]
\centering
    \begin{subfigure}{0.49\textwidth}
    	\centering
    	\includegraphics{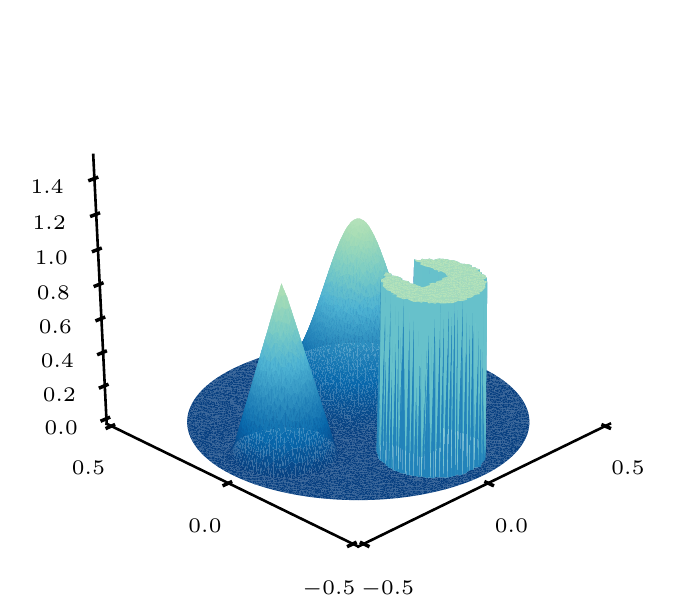}
    	\caption{Mesh 1, $t = 0$.}
    	%
    	%
    	\includegraphics{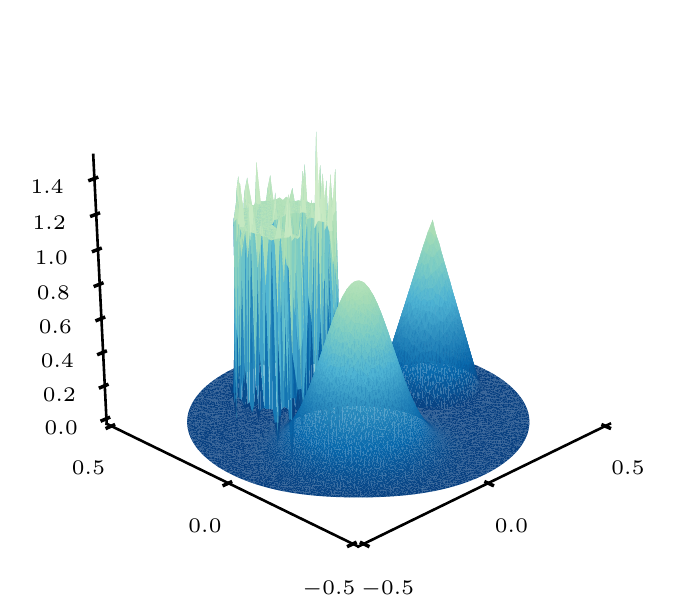}
    	\caption{Mesh 1, $t = 1$.}
    	%
    	%
    	\includegraphics{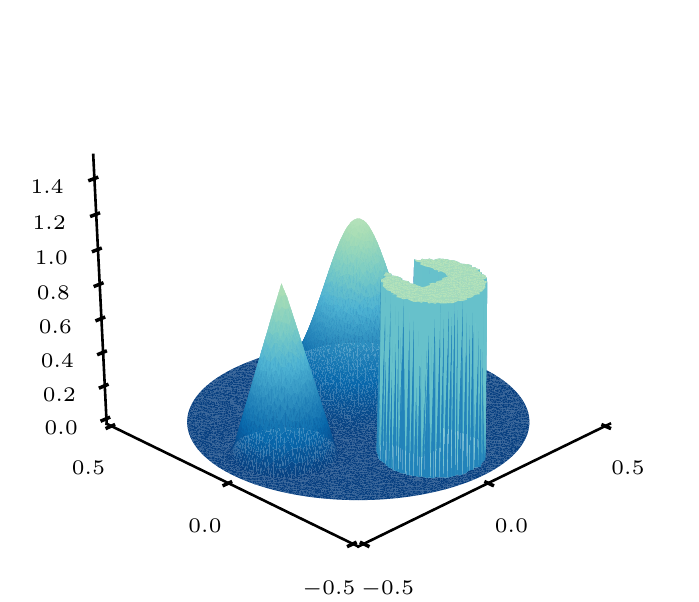}
    	\caption{Mesh 1, $t = 2$.}
    	\label{fig:mesh_1 - t=2}
    \end{subfigure}
    \begin{subfigure}{0.49\textwidth}
    	\centering
    	\includegraphics{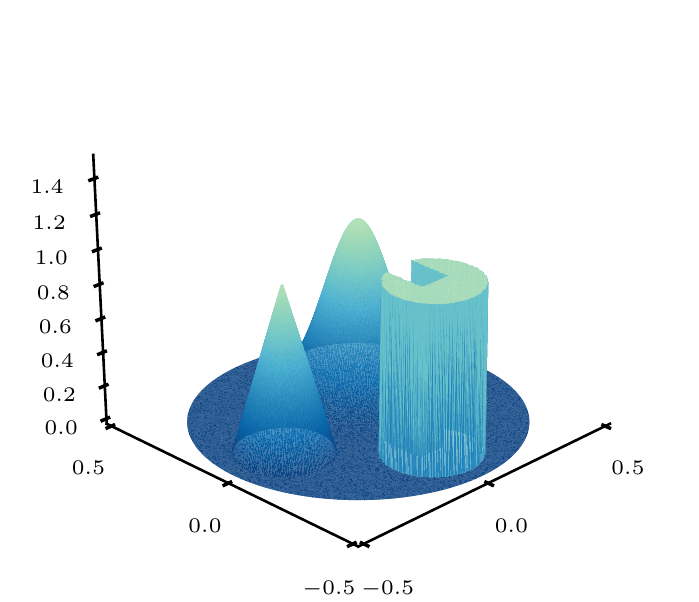}
    	\caption{Mesh 2, $t = 0$.}
    	%
    	%
    	\includegraphics{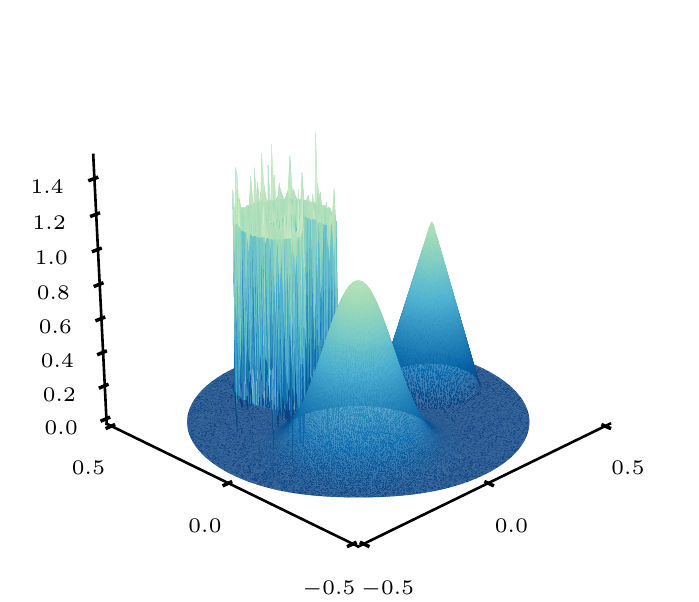}
    	\caption{Mesh 2, $t = 1$.}
    	%
    	%
    	\includegraphics{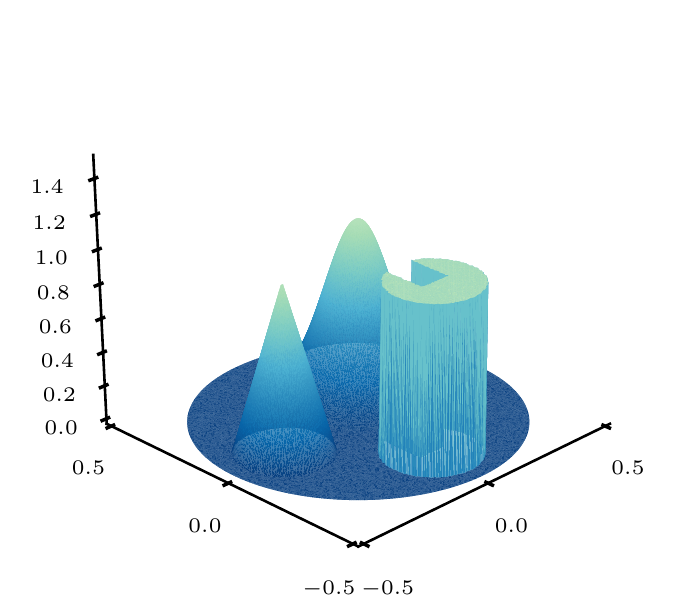}
    	\caption{Mesh 2, $t = 2$.}
    	\label{fig:mesh_2 - t=2}
    \end{subfigure}
    \caption{Rigid body rotation: numerical solution $\xDiscreteScalar{\psi}$ for different meshes at various time instants using polynomial orders $(k,l) = (1,0)$.}
    \label{fig:three_body_rotation}
\end{figure}

The results obtained for the two different configurations (see Table~\ref{tab:three_rigid_bodies}) are assessed visually after a half and a full rotation in Fig.~\ref{fig:three_body_rotation}. Since the particle values are not updated, the initial discontinuities are maintained at the particle level. As such, the test provides insight into how well the mesh-based solution follows the particle field. The shapes of the pointy cone and the Gaussian hump are well-preserved (without numerical damping), both for a half rotation as well as for a full rotation. 
Although the shape is well-preserved for the slotted-disk, localized overshoot is observed near the discontinuities for a half rotation. Rather than being a dispersion artifact, the over- and undershoot should be interpreted as a resolution issue with the mesh resolution being too coarse to capture the sharp discontinuity at the particle level monotonically. This is clearly illustrated by Fig.~\ref{fig:mesh_1 - t=2} and Fig.~\ref{fig:mesh_2 - t=2}, showing that the initial condition is accurately recovered for a full rotation at $t=2$, thus indicating that the mesh-based solution is able to follow the moving particle field.

We also investigate the mass conservation errors for the two different configurations. To this end, a measure for the relative global mass conservation error at time $T$ is defined as
    \begin{equation}\label{eq:area_error_gaussian}
    	\epsilon_{\Delta \phi_\Omega } =  \frac{ \areaIntegral{\Omega}{ \left( \phi_h( \mathbf{x}, T) - \phi_h(\mathbf{x}, 0) \right)} }{\areaIntegral{\Omega}{\phi_h(\mathbf{x}, 0)}},
    \end{equation} 
in which $\phi_h(\mathbf{x}, 0)$ and $\phi_h(\mathbf{x}, T)$ are the mesh related fields at time 0 and time $T$, respectively. \\
The local mass conservation error is investigated  via the $L^2$-norm of the time-discrete counterpart of the local conservation statement, Eq.~\eqref{eq:local_conservation}. For the problem under consideration, this local mass conservation error norm at time level $n+1$ is given by
    \begin{equation}\label{eq:area_error_local}
        \epsilon_{\Delta \phi_K} = \left( \sum_{K}^{} \left( \areaIntegral{K}{ \frac{\psi_h^{n+1} - \psi_h^{*,n}  }{\Delta t} } + \lineIntegral{ \partial K}{  \mathbf{a} \cdot \mathbf{n} \bar{\psi}_h^{n+1}  } \right)^2 \right)^{1/2}.
    \end{equation}
The mass conservation errors as defined by Eqs.~\eqref{eq:area_error_gaussian} and \eqref{eq:area_error_local} after a half rotation ($T=1$) and a full rotation ($T=2$) are tabulated in Table~\ref{tab:three_body_conservation}. The results confirm global and local mass conservation of the PDE-constrained projection to machine precision, irrespective of the resolution of the Eulerian mesh.
\begin{table}[H]
\centering
\caption{Rigid body rotation: mass conservation errors (defined by Eqs.~\eqref{eq:area_error_gaussian} and \eqref{eq:area_error_local}) for different mesh configurations after a half rotation ($T=1$) and a full rotation ($T=2$).}
\label{tab:three_body_conservation}
\small{
\begin{tabular}{c|c c| c c}
\hline
\hline 
				&   \multicolumn{2}{c|}{Mesh~1}  & \multicolumn{2}{c}{Mesh~2}\\
				&   T = 1   & T = 2 			& T = 1   & T = 2 \\
\hline 
$\epsilon_{\Delta \phi_\Omega}$  &  -2.0e-16 	& 2.0e-16 & -3.1e-15 & 5.9e-16	 \\
$\epsilon_{\Delta \phi_K}$       &  1.7e-16 	& 1.5e-16 & 1.3e-16        & 1.3e-16         \\
\hline 
\hline 
\end{tabular}
}
\end{table}

\subsection{Advection: advection skew to mesh}
Finally, we consider the advection of a discontinuity on the unit-square $\Omega := [0,1]^2$ for different transport velocities $\mathbf{u} = [\cos \alpha, \sin \alpha]^\top$ with characteristic directions $\alpha$ of $15 ^\circ, 30^\circ, 45^\circ, 60^\circ$. The diffusivity $\kappa$ is set to 0, so that we solve a pure advection problem. A regular triangular mesh is used with uniform cell size $h_K=1/25$, and each cell contains on average approximately 20 particles. Except for the case $\phi = 45 ^\circ$, the propagation directions are not aligned with the mesh. Dirichlet boundary conditions are prescribed at the inflow boundaries, and the specific mass $\psi_p$ carried by the particles flowing into the domain is set accordingly. Similar to the preceding rigid body rotation benchmark, the particle specific masses need not be updated, so that any discontinuities at the particle level persist without artificial diffusion. Hence, this test is well-suited for assessing the behavior of the scheme in the presence of steep gradients. Furthermore, in the particle-mesh setting of this study, the test case can be used to assess the global conservation statement in the presence of inflow and outflow boundary conditions. 
\begin{figure}[H]
	\centering
	\begin{subfigure}{0.49\textwidth}
		\centering
		\includegraphics[width = 0.85\textwidth]{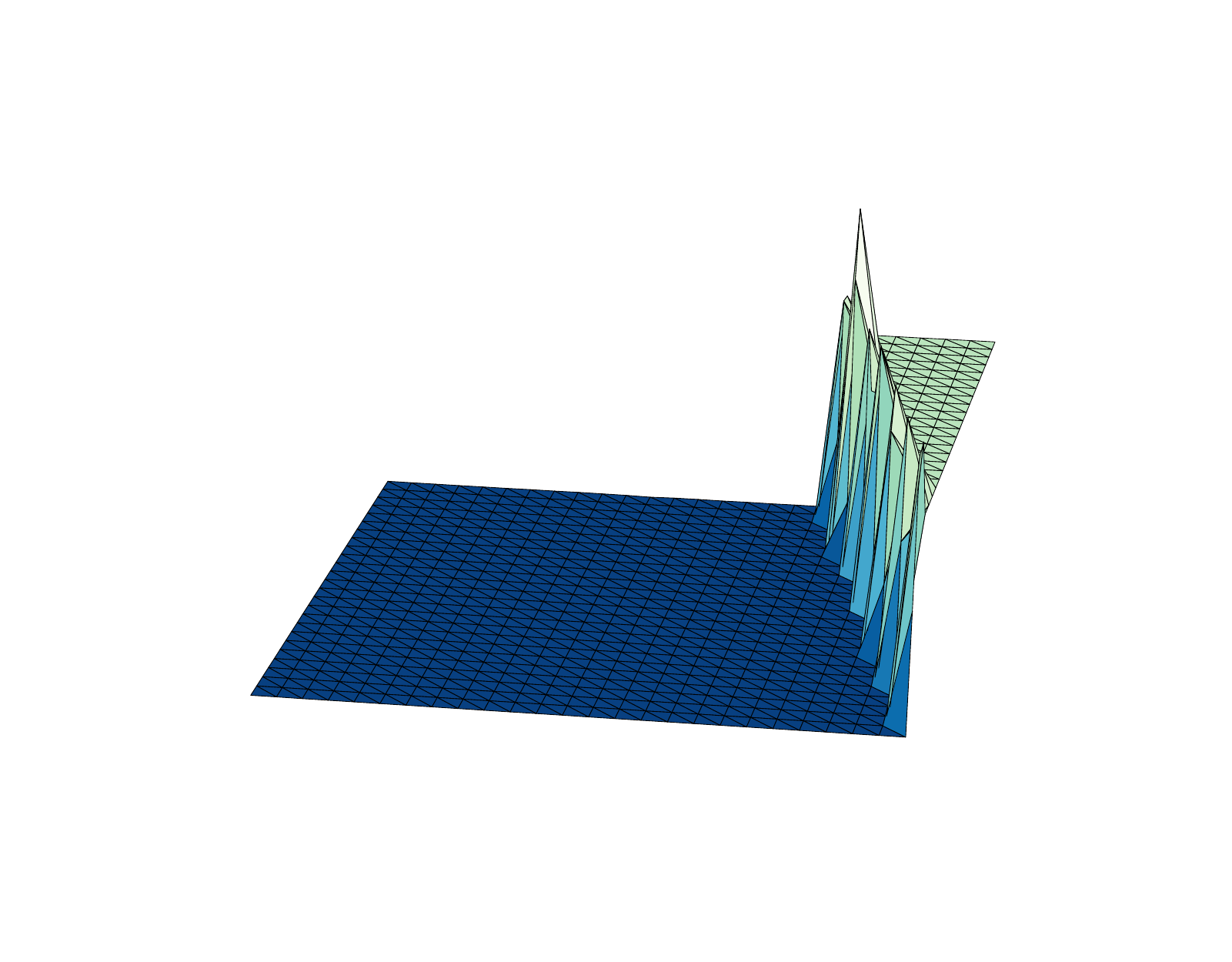}
		\caption{$\alpha = 15^\circ$.}
	\end{subfigure}
	\begin{subfigure}{0.49\textwidth}
		\centering 
		\includegraphics[width = 0.85\textwidth]{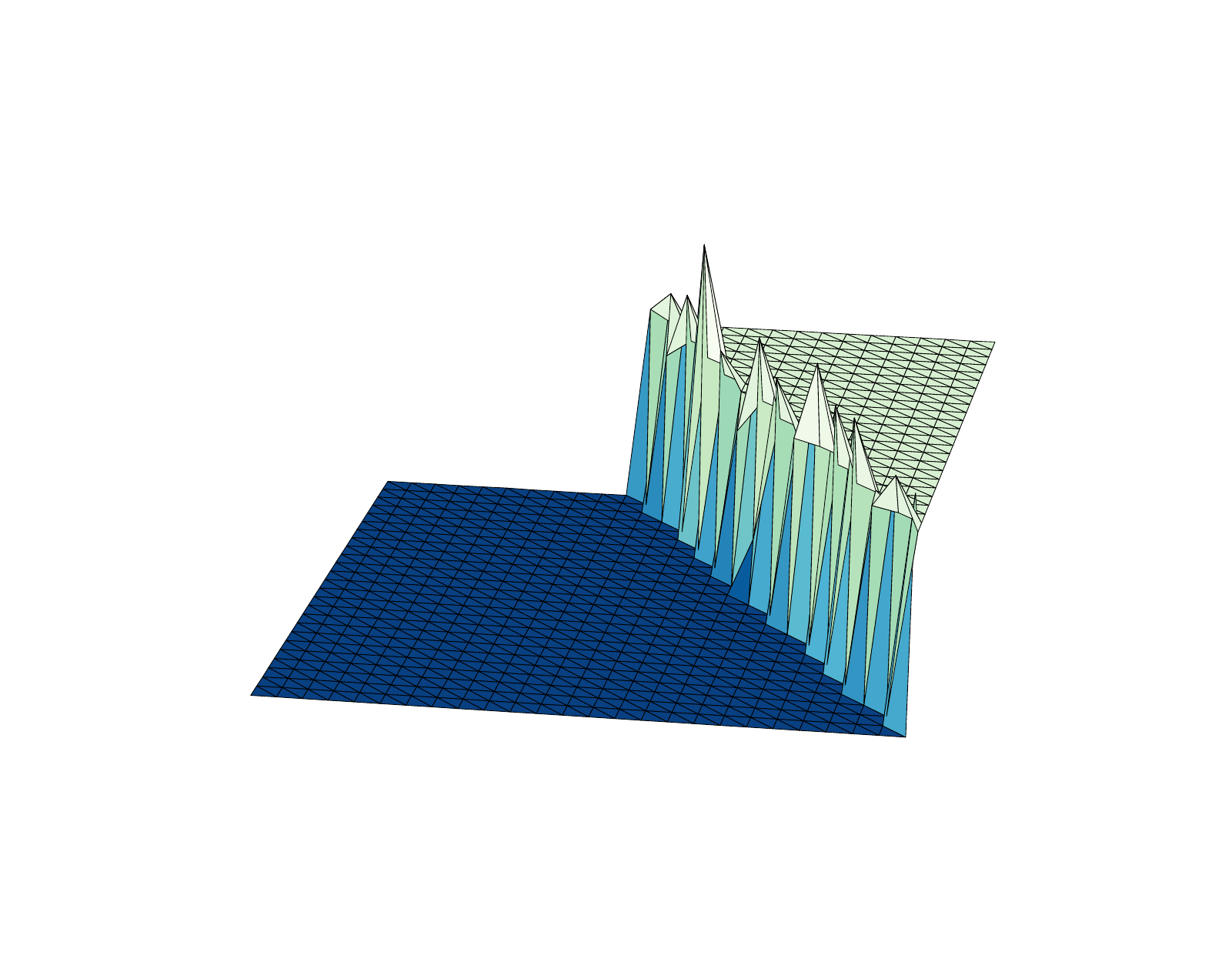}
		\caption{$\alpha = 30^\circ$}
		\label{fig:skew_advection_30}
	\end{subfigure}
	\begin{subfigure}{0.49\textwidth}
		\centering
		\includegraphics[width = 0.85\textwidth]{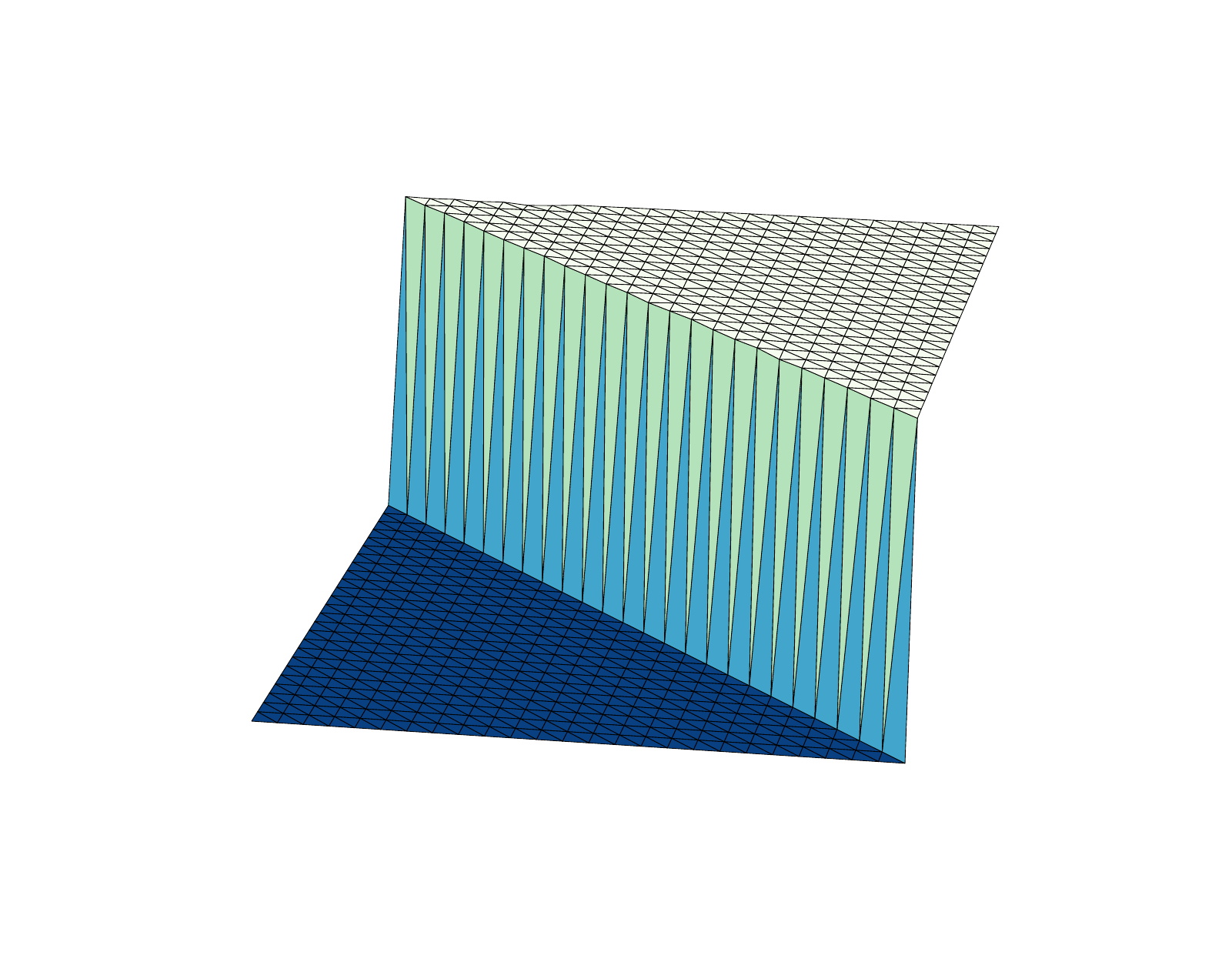}
		\caption{$\alpha = 45^\circ$}
	\end{subfigure}
	\begin{subfigure}{0.49\textwidth}
		\centering
		\includegraphics[width = 0.85\textwidth]{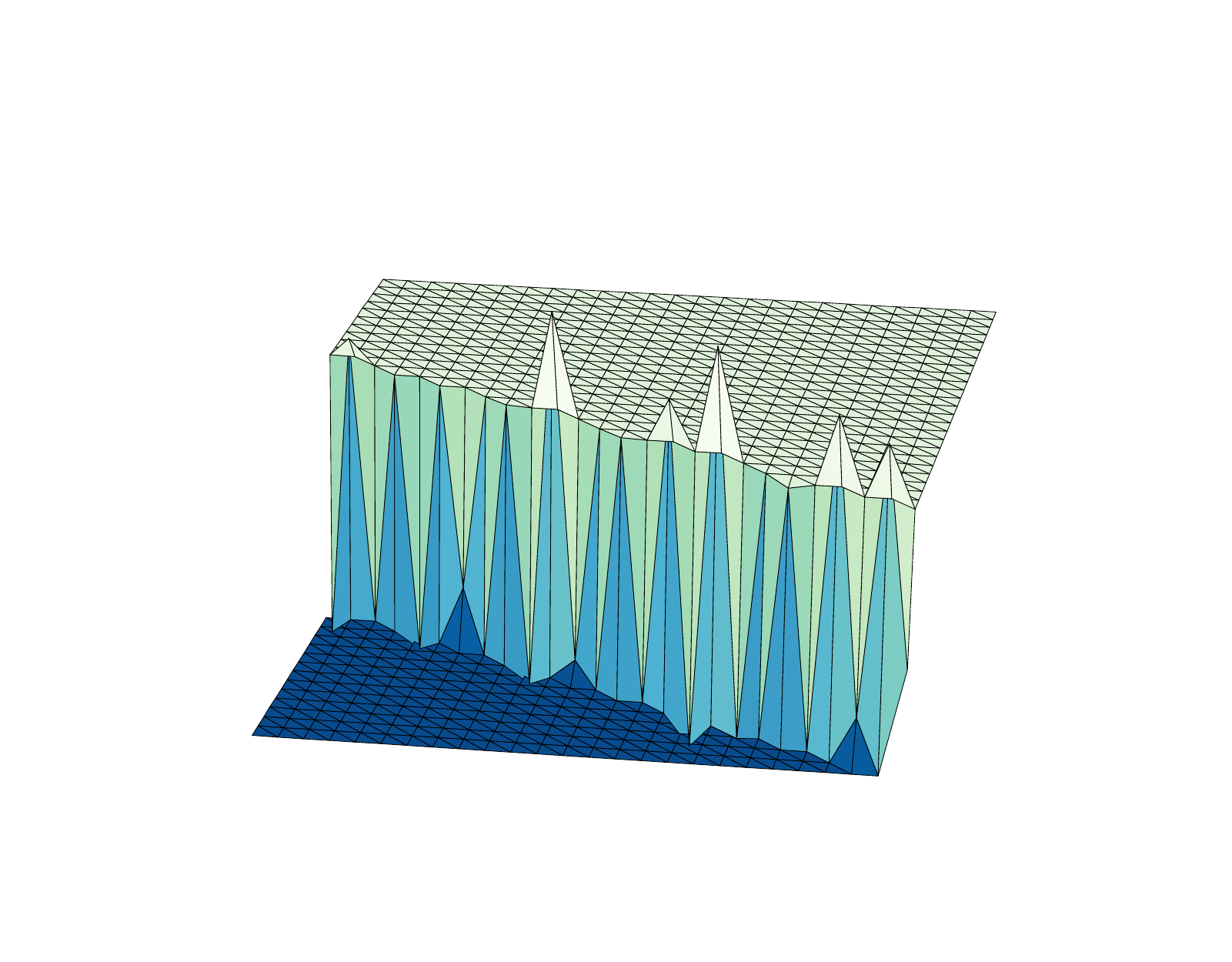}
		\caption{$\alpha = 60^\circ$}
	\end{subfigure}
	\caption{Advection skew to mesh: numerical solution $\xDiscreteScalar{\phi}$ on a unit square domain at $t=2.0$ for different characteristic directions using polynomial orders $(k,l)=(1,0)$ .}
	\label{fig:skew_to_mesh}
\end{figure}

The fields $\xDiscreteScalar{\phi}$ at $t=2.0$ are plotted in Fig.~\ref{fig:skew_to_mesh}. For a characteristic direction of $45^\circ$, the discontinuity is captured exactly at the mesh level. For the other characteristic directions, an overshoot is observed near the discontinuity. However, this overshoot remains strictly localized to one mesh cell upstream of the discontinuity. This behavior can be expected since no attempts are presently made to preserve monotonicity at the mesh, while the discontinuity at the particle level is inherently maintained without any diffusion. We leave the introduction of limiters as a fruitful area for future research, and refer to the work of Bochev and coworkers on bound-preserving remaps as an interesting starting point \cite{Bochev2013}. 

The global mass conservation property of the scheme is verified by virtue of the time-discrete equivalent of Eq.~\eqref{eq:global_conservation}. The global mass conservation error, denoted by $\epsilon_{\Delta \phi_\Omega}$, is the residual after subtracting the right-hand side from the left-hand side in this equation. The values for $\epsilon_{\Delta \phi_\Omega}$ thus obtained at $t=2$ are tabulated in Table~\ref{tab:skew_advection_conservation} for the four characteristic angles, confirming that global mass conservation is satisfied to machine precision.
\begin{table}[H]
\centering
\caption{Advection skew to mesh: global mass conservation errors for different characteristic directions.}
\label{tab:skew_advection_conservation}
\small{
\begin{tabular}{c|c|c|c|c}
\hline
\hline 
$\alpha$					&   $15^\circ$  & $30^\circ$ & $45^\circ$ & $60^\circ$ \\
\hline 
$\epsilon_{\Delta \phi_\Omega}$  &  -1.25e-15 	& -2.70e-15	 & -4.43e-15  & -4.75e-15 \\
\hline 
\hline 
\end{tabular}
}
\end{table}

\subsection{Navier-Stokes equations: Poiseuille flow}
Plane Poiseuille flow is considered to assess the behavior of the particle-mesh operator splitting scheme for the incompressible Navier-Stokes equations in the presence of no-slip boundaries. In addition, the test is used to elaborate upon the required minimum number of particles per cell. 

Starting from rest, the flow gradually develops towards a steady state under the influence of a constant, axially applied body force $ \mathbf{f} $. The channel is modeled in the $(x,y)$-plane, with the $x$-axis pointing in the flow direction. The model domain is given by $\Omega := [0, 1] \times [-0.25, 0.25]$. The flow is periodic in the $x$-direction and at the location of the plates ($y=\pm d$) no-flux and no-slip boundary conditions are used in the PDE-constrained particle-mesh projection and the Stokes step, respectively. The kinematic viscosity is set to $\nu = 1\cdot10^{-3}$, and the body force $\mathbf{f} = (f_x,0)^T$ is chosen such that the steady state Reynolds number $Re = 2Ud/\nu$ equals 200. The analytical transient solution of the axial velocity can be found in, e.g., \cite{Sigalotti2003}. 
\subsubsection{Convergence study}
A coarse, regular triangular mesh is constructed using $8\times4\times2 = 64$ cells. The mesh size is halved three times, giving a sequence of mesh refinements with the finest mesh containing 4096 cells. On average, 30 particles are initially assigned per cell. This number is kept constant upon mesh refinement. Throughout, we will use a $\theta_L$-value of 1/2 for the particle updating, and the time step size corresponds to a CFL-number of approximately 0.25. Results at the dimensionless time instant $t^* = t U/2d = 100$ are presented in Table \ref{tab:posseuille_convergence_bc} for polynomial orders $(k,l)=(1,0)$ and $(k,l)=(2,0)$, respectively. The observed convergence is near-optimal (i.e. order $k+1$ in the velocity and order $k$ in the pressure).
\begin{table}[H]
\centering
\caption{Poiseuille flow: convergence of the $L^2$-error in the velocity and the pressure at dimensionless time $t*= t U/2d = 100$ for different polynomial orders $(k,l)$.} 
\label{tab:posseuille_convergence_bc}
\small{
\begin{tabular}{c c c|c c c c}
\hline
\hline
$(k,l)$ & Cells & $ \Delta t $ & $\lVert \mathbf{u} - \mathbf{u}_h \rVert$ &  Rate &  $\lVert p - p_h \rVert$&  Rate \\
\hline
\multirow{ 4}{*}{$(1,0)$} 
&     64 &      0.2  & 6.2e-3  & -		&	1.4e-4 	&	-	\\ 
&    256 &      0.1  & 1.6e-3  & 2.0	&	8.1e-5	&	0.7 \\
&   1024 &     0.05  & 3.9e-4  & 2.0	&	3.7e-5	&	1.2 \\
&   4096 &    0.025  & 9.7e-5  & 2.0	&	1.8e-5	&	1.0 \\
\hline
\multirow{ 4}{*}{$(2,0)$} 
&     64 &       0.2 & 3.9e-6  &  -   & 1.3e-7  & -   \\ 
&    256 &      0.1  & 4.3e-7  &  3.2 & 1.8e-8  & 2.8 \\ 
&   1024 &     0.05  & 5.1e-8  &  3.1 & 3.7e-9  & 2.3 \\ 
&   4096 &    0.025  & 5.1e-9  &  3.3 & 8.0e-10 & 2.2 \\ 
\hline 
\hline
\end{tabular}
}
\end{table}
\subsubsection{Assessing the particle resolution}
Evidently, it is desirable from an efficiency perspective to keep the number of particles as low as possible without compromising accuracy. We therefore investigate the influence of the particle resolution on the accuracy by considering the Poiseuille flow benchmark, using $l=0$ and $k=1,2,3,4$, combined with a variable particle resolution. We restrict the discussion to the mesh containing 256 cells and a time step of $\Delta t = 0.1$ and only vary the particle resolution for this configuration so that the number of particles per cell (denoted with $\Bar{S}_0^K$) is in the range 2-50, initially. In order to have full control over the initial particle configuration, particles are placed on a regular lattice. 
\begin{figure}[H]
\centering
\includegraphics{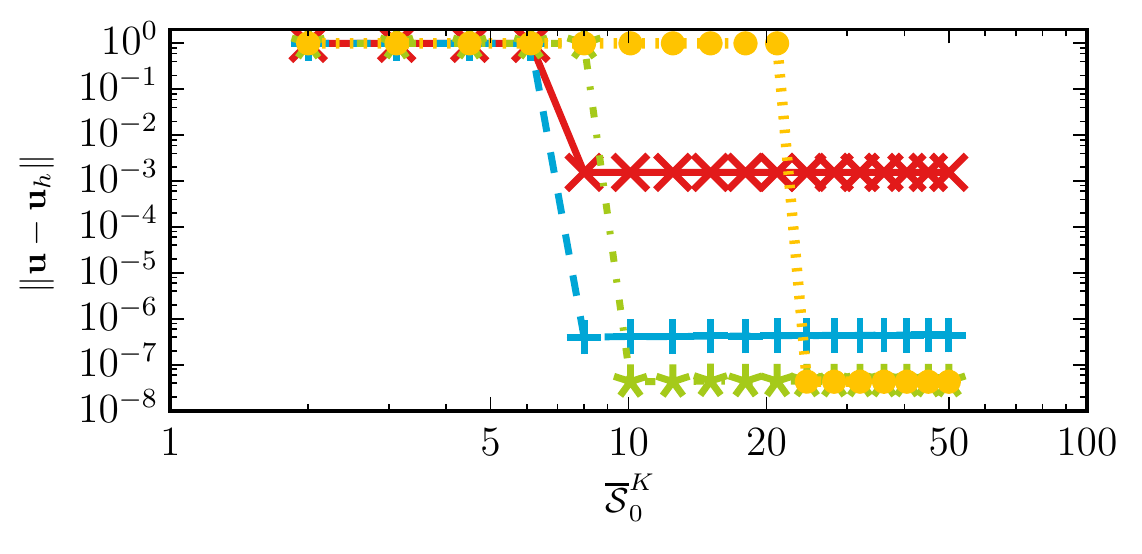}
\caption{Poiseuille flow: $L^2$-error at $t^* = 100$ as a function of the average number of particles initially assigned per cell ($\overline{\mathcal{S}}_0^K$) for an Eulerian mesh with a fixed number of 256 cells and a time step of $\Delta t = 0.1$ for a linear ($\times$), quadratic ($+$), cubic ($*$) and quartic ($\bullet$) polynomial basis.}
\label{fig:particle_resolution test}
\end{figure}
As expected, many of the low particle resolution tests fail prematurely (i.e. before reaching the end time $t^* = t U/2d = 100$). In Fig.~\ref{fig:particle_resolution test}, results are visualized by plotting the $L^2$-error of the velocity field against the (average) number of particles per cell for the different polynomial orders. For convenience, we assign an error-value of 1 to all low particle resolution runs failing prematurely.

Evident from the figure is the sharp transition in the error levels. That is, for low particle resolutions the computations are prone to fail prematurely, whereas from a certain particle resolution threshold onward, accurate results are obtained with error values being independent of the particle resolution. Interestingly, this threshold particle resolution depends on the polynomial order of the basis functions and is approximately equal to $\Bar{S}_0^K=8, 8, 10, 15$ for $k=1,2,3,4$, respectively. 
By recognizing that the particle-mesh projection is based on the local least-squares minimization problem (Eq.~\eqref{eq:J1}), we note that unisolvency of the particle locations with respect to the discontinuous function space $\xDiscreteVector{W}$ is a necessary condition for the particle-mesh projection to be accurate. Empirically, we thus observe this condition to be satisfied if the minimal number of particles in a cell is at least equal to the number of local basis functions, except for linear basis functions ($k=1$) and to a lesser extent for quadratic basis functions ($k=2$), which require a higher number of particles per cell. 
\subsection{Navier-Stokes equations: Taylor-Green flow}
The Taylor-Green flow on the bi-periodic domain $\Omega := [-1, 1] \times [-1, 1]$ is considered as a final example. This problem features a periodic sequence of decaying vortices. Provided the Reynolds number is sufficiently small, closed analytical expressions for the velocity and the pressure are given by, respectively,
	\begin{align} 
	\mathbf{u}(\mathbf{x},t) &= U\exp\left(-2\nu\pi^2t\right)\left(-\cos(k_x x)\sin(k_y y)\hspace{3pt}, \hspace{3pt} \sin(k_x x)\cos(k_y y) \right)^\top , \label{eq:tg_analytic_velocity} \\ 
	p(\mathbf{x},t) &= \frac{1}{4}\exp\left(-4\nu\pi^2t\right)\left(\cos(2k_x x)+\cos(2k_y y)\right) ,\label{eq:tg_analytic_pressure}
	\end{align}
in which $U$ is the initial velocity amplitude, $\nu$ is the kinematic viscosity, and $k_x=2\pi/L_x$ and $k_y=2\pi/L_y$ are wave numbers in the $x$- and $y$-direction, respectively, with $L_x$ and $L_y$ being the associated wave lengths.\\
Since only an initial condition needs to be specified, this test is particularly suited to study the accuracy of the time stepping scheme. Furthermore, due to the absence of body forces and boundary tractions the test is useful to investigate momentum conservation. To this end, a mesh-related measure for the global momentum conservation error is defined as 
	\begin{align}
	\epsilon_m &= \left| \areaIntegral{\Omega}{ \left(\mathbf{u}_h(\mathbf{x}, t)-\mathbf{u}_h(\mathbf{x}, 0) \right)} \right|. \label{eq:global_momentum_error}
	\end{align} 
\subsubsection{Convergence study}
To investigate the convergence properties of the scheme, we consider a time interval of interest $I=(0,2]$. The initial peak velocity $U$  in Eq.~\eqref{eq:tg_analytic_velocity} is set to 1, the wave lengths $L_x = L_y = 2$, and the kinematic viscosity is either $\nu = 2\cdot 10^{-2}$ or $\nu = 2\cdot 10^{-3}$, resulting in Reynolds numbers of $Re = UL/\nu = 100$ and $Re = 1000$, respectively. The chosen time step corresponds to a CFL-number of approximately 0.4. Different function space combinations are considered, see Table \ref{tab:taylor_green_overview}. For comparison, Case~3 uses the (non-conservative) $\ell^2$ particle-mesh interaction from \cite{Maljaars2017}.
\begin{table}[H]
\centering
\caption{Taylor-Green flow: Overview of model settings.}
\label{tab:taylor_green_overview}
\small{
\begin{tabular}{c|c c c c}
\hline
\hline 
\rule{0pt}{2.5ex}   
		& Projection Method	&$k$	&	$l$	&	$\theta_L$	\\
\hline 
Case 1 	& PDE		& 1	& 	0	& 		1/2		\\
Case 2 	& PDE		& 2	& 	0	& 		1/2		\\
Case 3 	& $\ell^2$	& 2	& 	0	& 		1/2		\\
\hline 
\hline 
\end{tabular}
}
\end{table}
\begin{table}[H]
\centering
\caption{Taylor-Green flow: overview of model runs with the associated errors  $\lVert \mathbf{u} - \mathbf{u}_h \rVert$, $\lVert p - p_h \rVert$ and $\epsilon_{m}$ at time $t=2$.} 
\label{tab:taylor_green_errors_theta_imex}
\small{
\begin{tabular}{c r r|c c c c c|c c c c c}
\hline
\hline
\multicolumn{3}{c|}{} &  \multicolumn{5}{c|}{$Re = 100$}  & \multicolumn{5}{c}{$Re = 1000$}\\
\hline
& \scriptsize{Cells}  & \scriptsize{$\Delta t$}  & \scriptsize{$\lVert \mathbf{u} - \mathbf{u}_h \rVert$} &  \scriptsize{Rate} & \scriptsize{$\lVert p - p_h \rVert$} &  \scriptsize{Rate} & $\epsilon_m$  & \scriptsize{$\lVert \mathbf{u} - \mathbf{u}_h \rVert$} &  \scriptsize{Rate} &  \scriptsize{$\lVert p - p_h \rVert$}  &  \scriptsize{Rate} & $\epsilon_m$ \\
\hline
\multirow{ 4}{*}{Case 1}
		&    128 &	0.1 	&	1.2e-1 & - 		& 6.5e-2 & - 	& 2.8e-15 & 2.4e-1 & - 		& 2.9e-1 & - 	& 1.5e-15\\ 
		&    512 &	0.05	& 	2.4e-2 & 2.3 	& 2.6e-2 & 1.3 	& 4.5e-15 & 4.8e-2 & 2.3 	& 8.7e-2 & 1.7 	& 8.1e-15\\ 
		&   2048 &	0.025 	& 	4.7e-3 & 2.4	& 1.2e-2 & 1.1 	& 3.5e-15 & 1.1e-2 & 2.1 	& 4.0e-2 & 1.1 	& 3.6e-15\\ 
		&   8192 &	0.0125 	& 	1.5e-3 & 1.6 	& 6.0e-3 & 1.0 	& 1.2e-14 & 4.0e-3 & 1.5 	& 2.0e-2 & 1.0 	& 7.1e-16\\ 
\hline
\multirow{ 4}{*}{Case 2}
		&    128 &	0.1 	& 6.5e-3 & - 	& 1.5e-2 & -	& 8.8e-14 & 7.6e-2 & - 		& 5.6e-2 & - 	& 1.4e-13 \\  
		&    512 &	0.05	& 1.9e-3 & 1.8 	& 3.2e-3 & 2.2	& 1.6e-13 & 1.2e-2 & 2.7 	& 1.3e-2 & 2.1 	& 3.1e-13 \\ 
		&   2048 &	0.025   & 5.1e-4 & 1.9 	& 8.4e-3 & 1.9	& 3.4e-13 & 2.3e-3 & 2.4 	& 3.2e-3 & 2.0 	& 6.1e-13 \\ 	
		&   8192 &	0.0125  & 1.3e-4 & 2.0 	& 2.1e-4 & 2.0	& 6.3e-13 & 5.6e-4 & 2.0 	& 7.8e-4 & 2.0 	& 1.3e-12 \\   
\hline
\multirow{ 4}{*}{Case 3}
		&    128 &	0.1 	& 6.6e-3 & -    & 1.5e-02 & -   & 3.0e-4 & 7.6e-2 & -  		& 5.6e-2 & -  & 1.3e-3 \\
		&    512 &	0.05	& 1.9e-3 & 1.8 	& 3.2e-03 & 2.2 & 1.2e-4 & 1.2e-2 & 2.7 	& 1.3e-2 & 2.1 & 6.4e-5 \\
		&   2048 &	0.025 	& 5.2e-4 & 1.9 	& 8.5e-04 & 1.9 & 2.9e-6 & 2.2e-3 & 2.4 	& 3.1e-3 & 2.1 & 1.4e-6 \\ 	
		&   8192 &	0.0125	& 1.3e-4 & 2.0 	& 2.1e-04 & 2.0 & 2.8e-7 & 5.5e-4 & 2.0 	& 7.6e-4 & 2.0 & 1.5e-6 \\       
\hline
\hline
\end{tabular} 
}
\end{table}

Velocity and pressure errors are tabulated in Table~\ref{tab:taylor_green_errors_theta_imex} for different model runs. In all these runs, the average number of particles per cell is 28. For Case 1, consistent second-order convergence is observed for the velocity and first-order convergence is obtained for the pressure. Given the function spaces (piecewise linear velocity, and piecewise constant pressure), no better convergence rates would be expected. \\
For Case 2, we obtain approximately second-order convergence in the velocity and the pressure, both for the $Re =100$ test case and the $Re = 1000$ test case. 
This indicates that the time stepping error becomes dominant over the spatial error, where the former is expected to converge with second-order. 
As expected, momentum is conserved globally up to machine precision for the PDE-constrained particle-mesh interaction (Case 1, Case 2), whereas this is clearly not so for the unconstrained projection method (Case 3), using local $\ell^2$-projections. Noteworthy to mention is that the errors in the velocity and pressure, respectively, are almost identical for Case 2 and Case 3. 
\subsubsection{High Reynolds number test} 
\label{sec:taylor_green_low_reynolds}
Finally, we show that the developed scheme is robust for high Reynolds numbers. To this end, we investigate an $8\times8$ array of Taylor-Green vortices (i.e. using $L_x=L_y =0.5$). 
The initial velocity amplitude is $U= 1$ and the kinematic viscosity is chosen such that the Reynolds number is $Re = U L /\nu = 2.5\cdot10^5$. A $200\times200\times2$ regular triangular mesh is used, the polynomial order are set to $k=2$ and $l=0$, and $2.25\cdot 10^6$ particles are seeded in the domain resulting in approximately 28 particles per cell, initially. Furthermore, the time step size is $\Delta t = 5\cdot10^{-3}$, resulting in a CFL-number of approximately 1.
\begin{figure}[H]
	\centering
	\begin{subfigure}{0.3\textwidth}
		\centering
		\includegraphics[width = \textwidth]{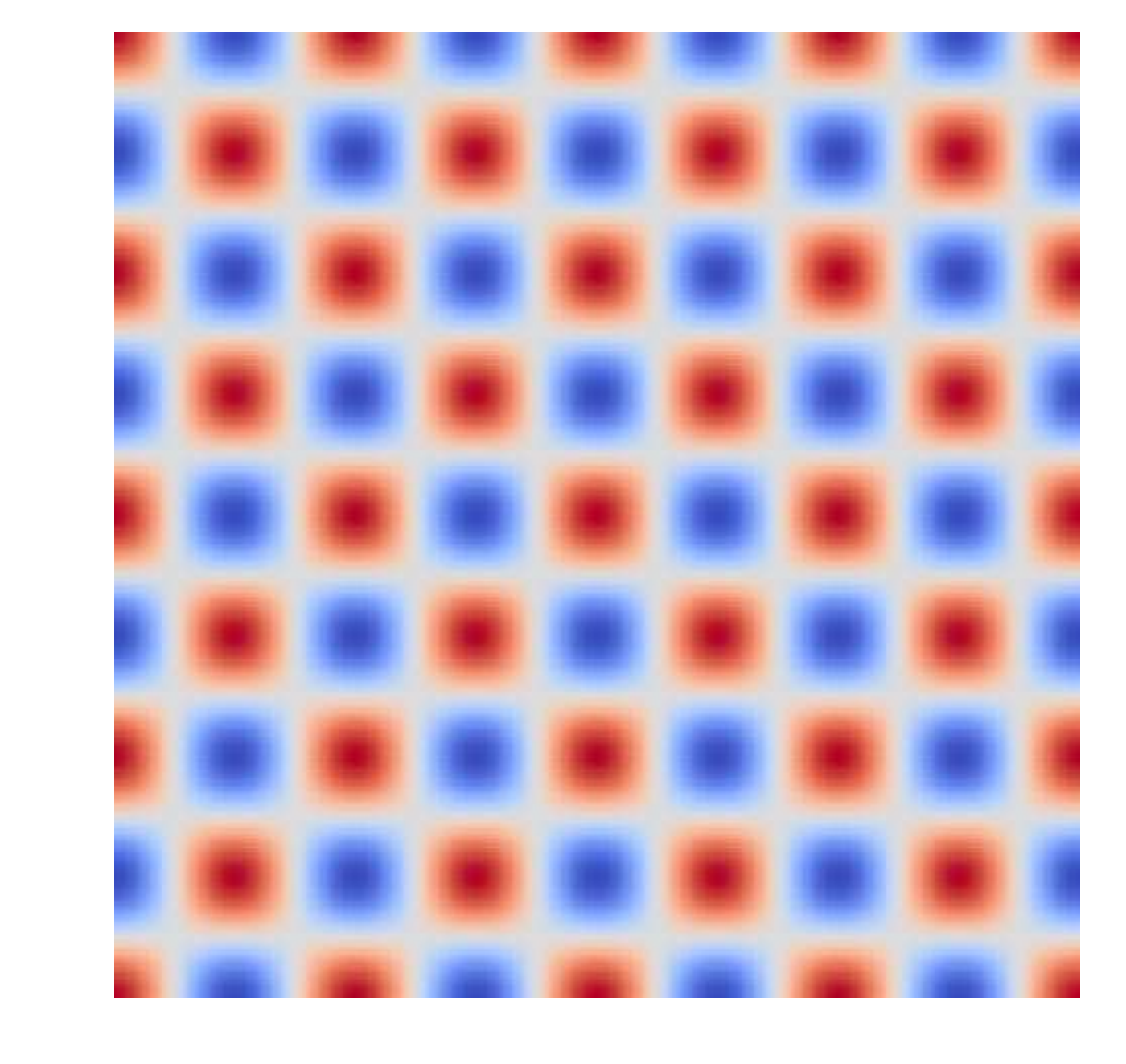}
		\caption{$t=0.0$}
	\end{subfigure}
	\begin{subfigure}{0.3\textwidth}
		\centering
		\includegraphics[width = \textwidth]{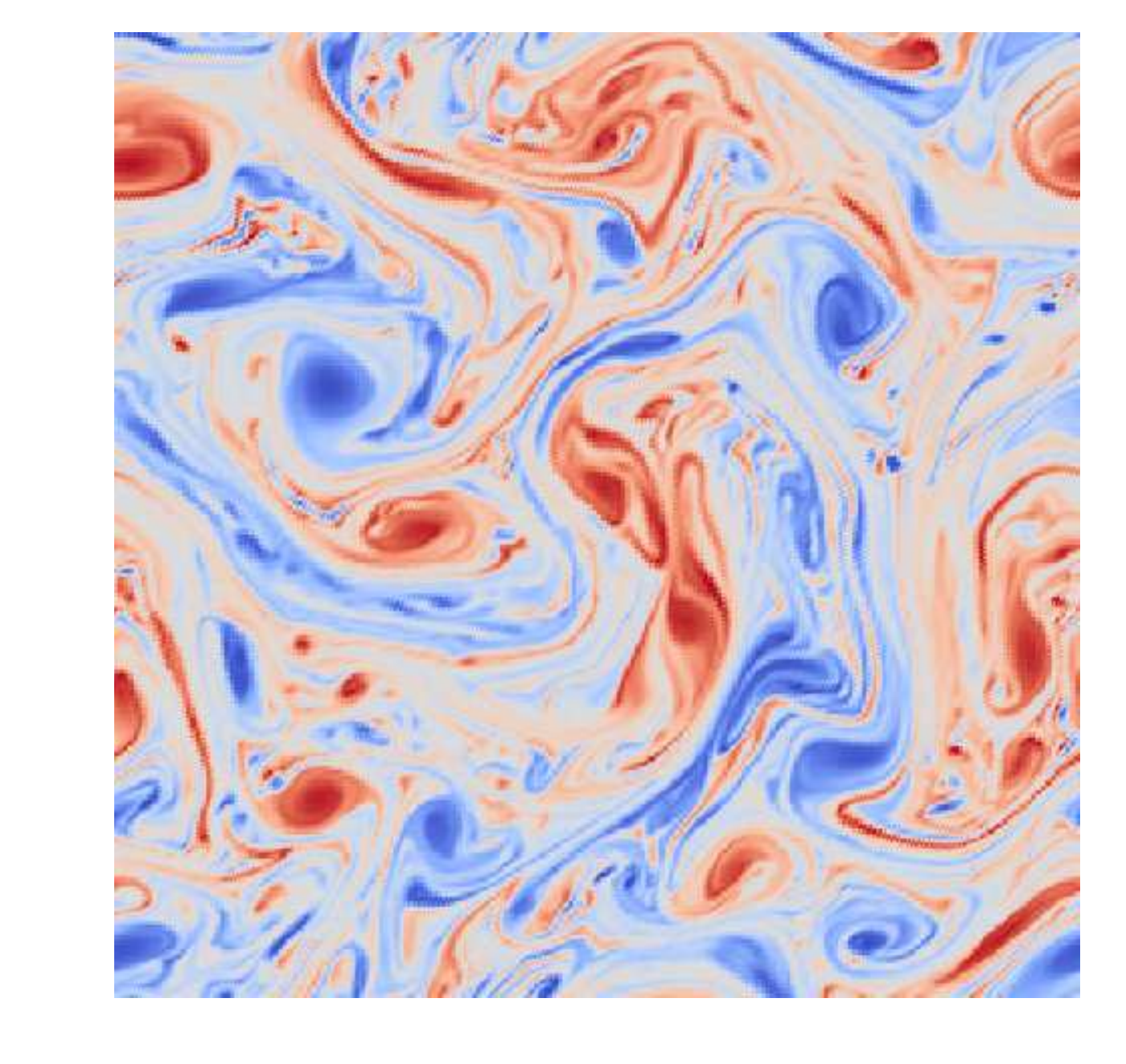}
		\caption{$t=7.5$}
	\end{subfigure}
	\begin{subfigure}{0.3\textwidth}
		\centering 
		\includegraphics[width = \textwidth]{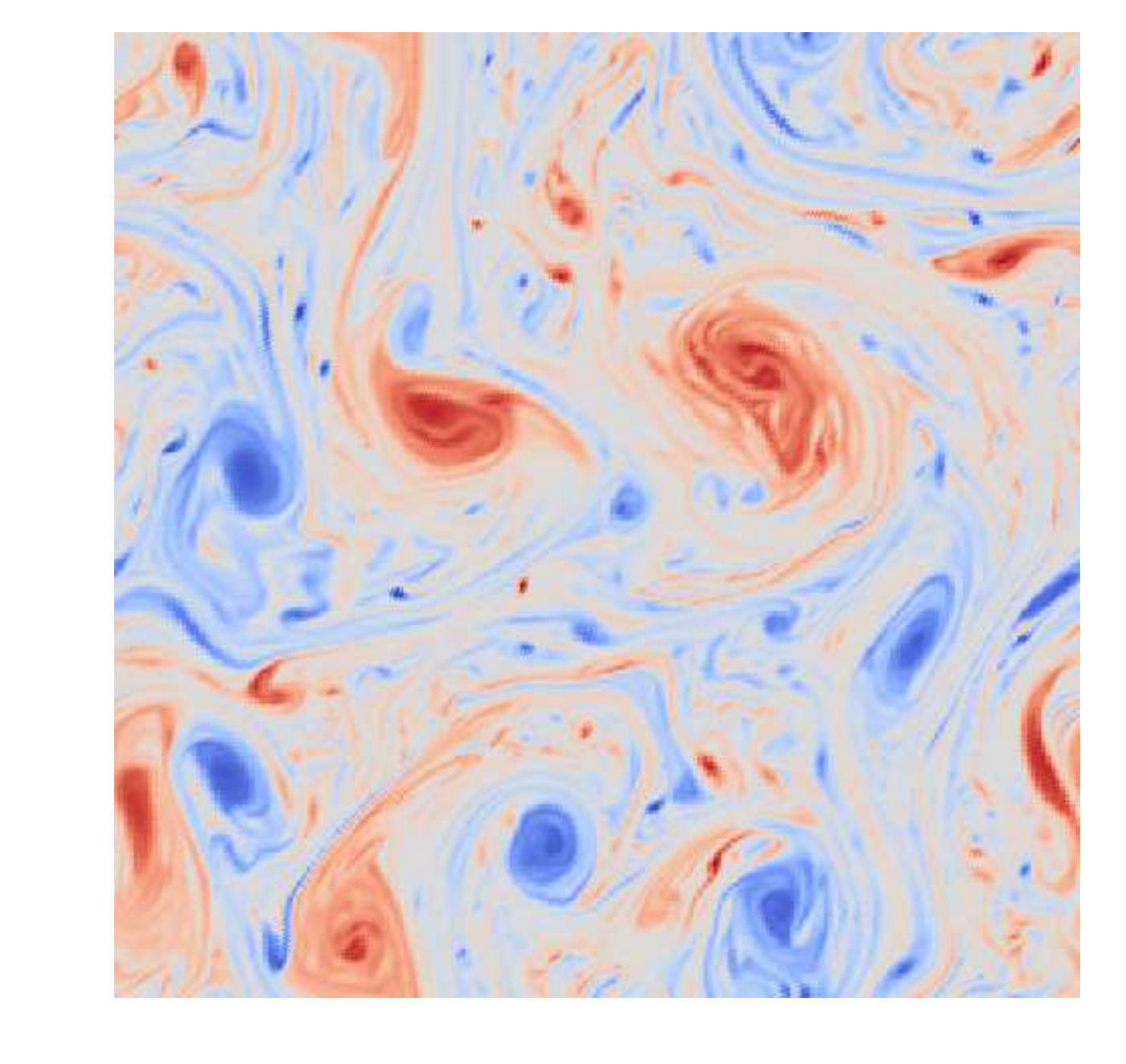}
		\caption{$t=10.0$}
	\end{subfigure} \\
	\begin{subfigure}{0.3\textwidth}
		\centering
		\includegraphics[width = \textwidth]{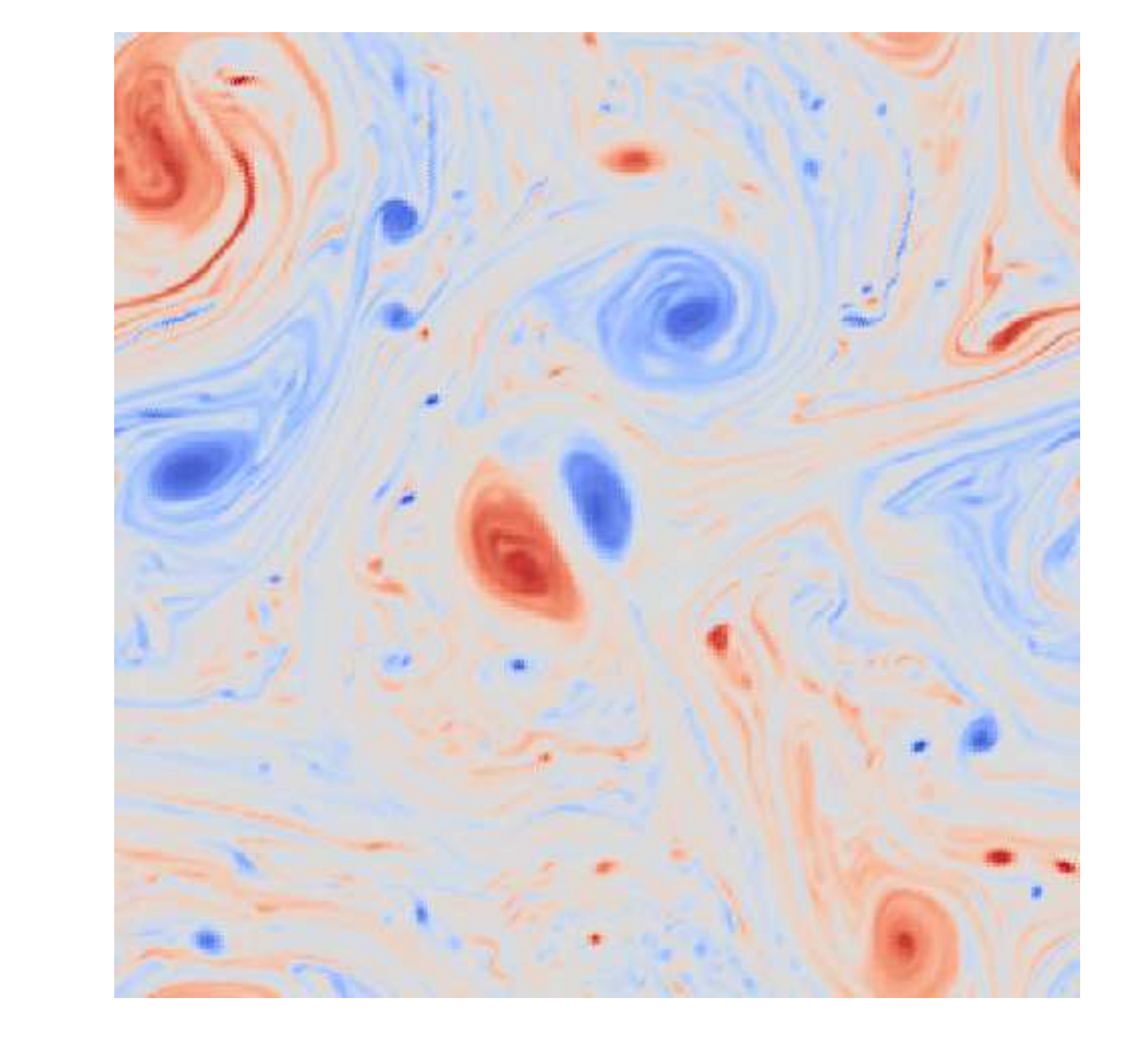}
		\caption{$t=12.5$}
	\end{subfigure}
	\begin{subfigure}{0.3\textwidth}
		\centering 
		\includegraphics[width = \textwidth]{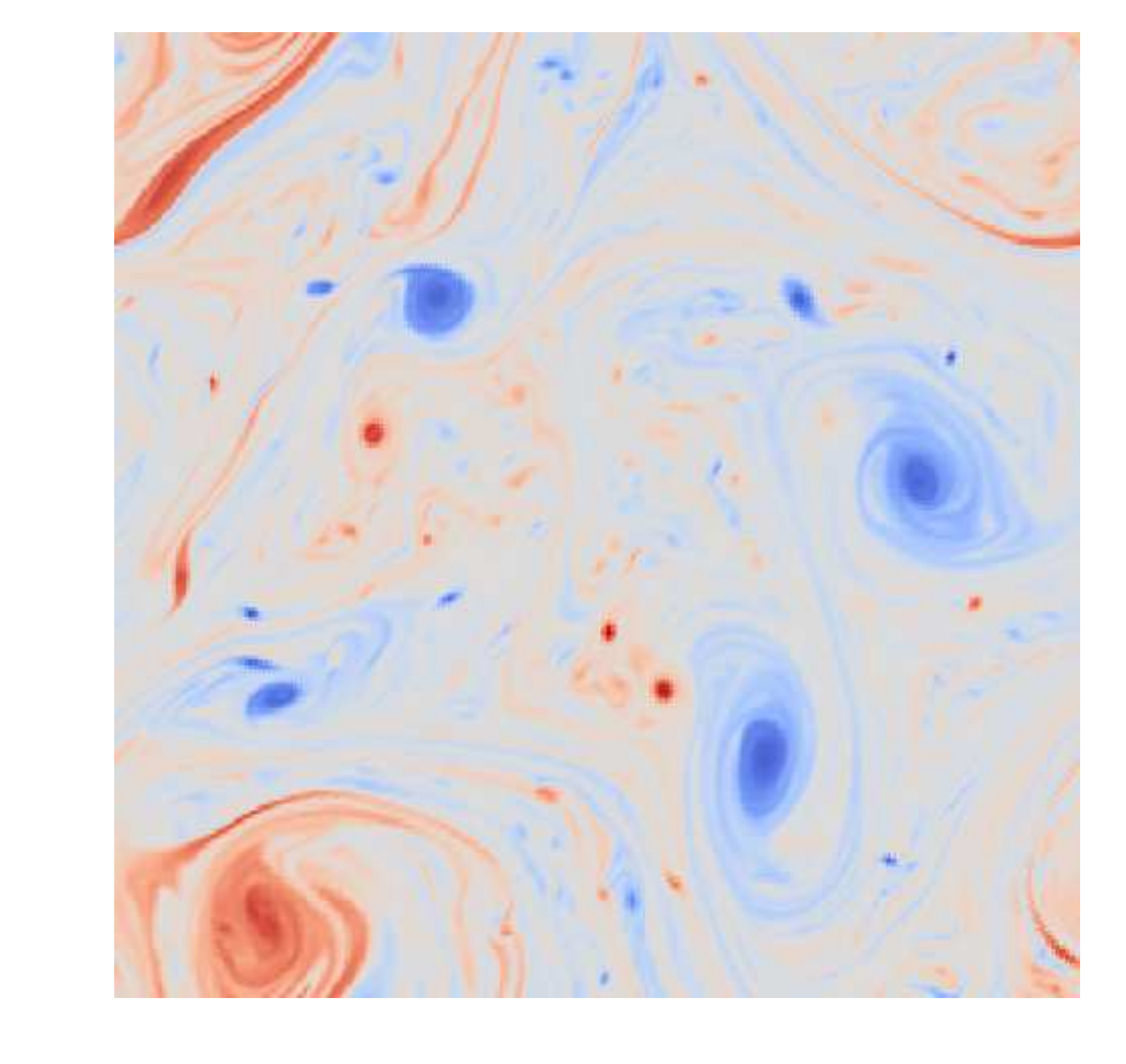}
		\caption{$t=15.0$}
	\end{subfigure}
		\begin{subfigure}{0.3\textwidth}
			\centering 
			\includegraphics[width = \textwidth]{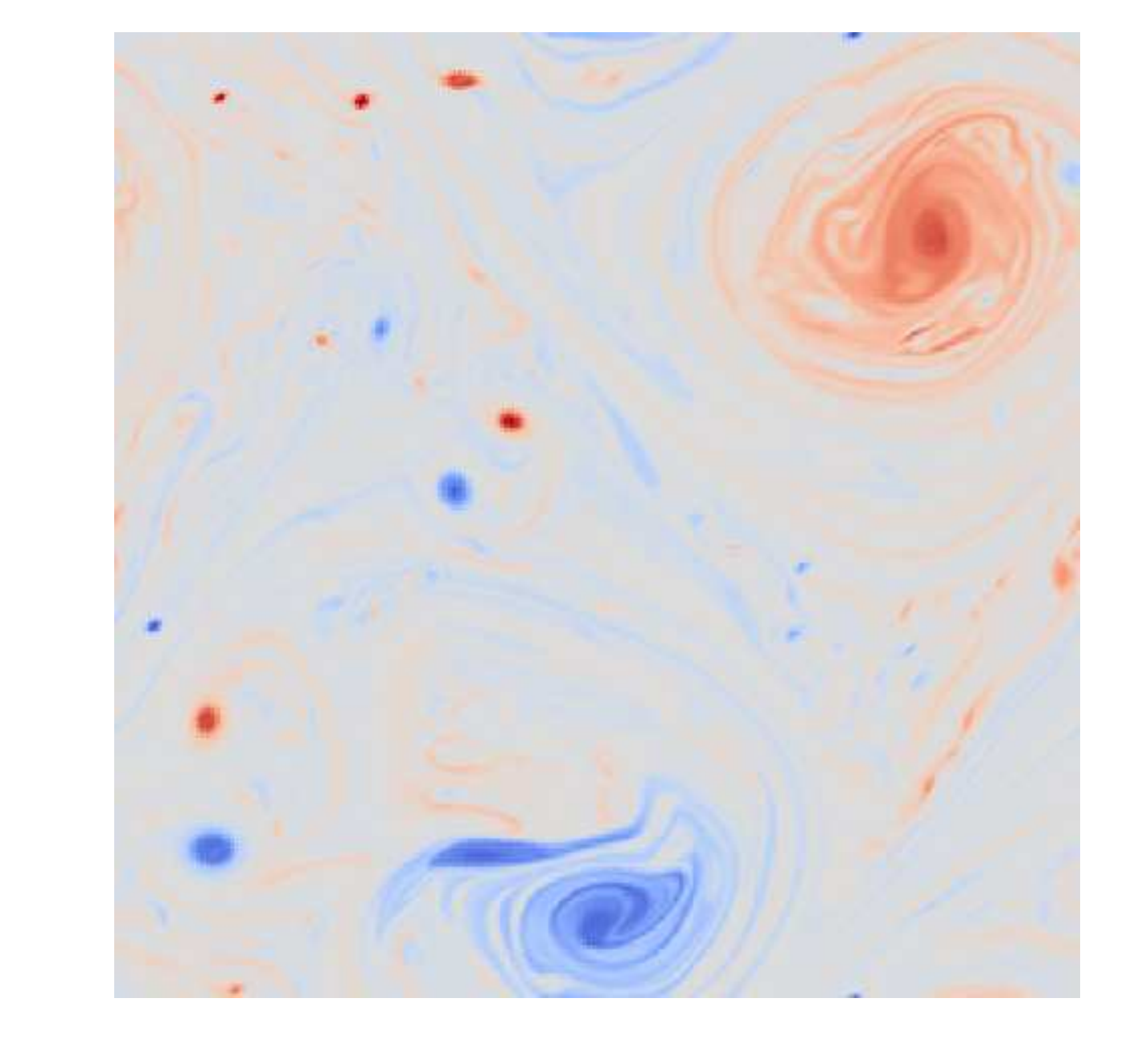}
			\caption{$t=20.0$}
		\end{subfigure}
	\caption{Taylor-Green flow: Vorticity fields for Taylor-Green flow with $Re=2.5\cdot10^5$ at different time instances. Color scaling is between -25 (blue) and 25 (red).}
	\label{fig:taylor_green_lonu}
\end{figure}

The vorticity field at various time instants is plotted in Fig.~\ref{fig:taylor_green_lonu}. For this high Reynolds number, the regular initial lattice of vortices collapses into a chaotic, turbulent motion. Two features characteristic of 2D turbulence are qualitatively observed in Fig.~\ref{fig:taylor_green_lonu}. First, comparing the vorticity fields at $t=7.5$ and $t=20$, the emergence of large-scale coherent structures is clearly observed for longer simulation runtimes. These large-scale vortices are only weakly dissipative and hence persist for many eddy-turnover times \cite{Clercx2009}. Second, the viscous dissipation takes place in the thin vortex filaments, clearly observed in Fig.~\ref{fig:taylor_green_lonu}.
Somewhat remarkable are the small-scale, high-vorticity spots, best observed for the longer simulation runtimes (e.g. for $t=15$ or $t=20$). As reported in literature \cite{Tabeling2002}, such small-scale vortices can be generated by Kelvin-Helmholtz type instabilities of the vortex filaments for sufficiently high Reynolds numbers.

\begin{figure}[H]
	\centering
	\begin{subfigure}{0.495\textwidth}
		\centering
		\includegraphics{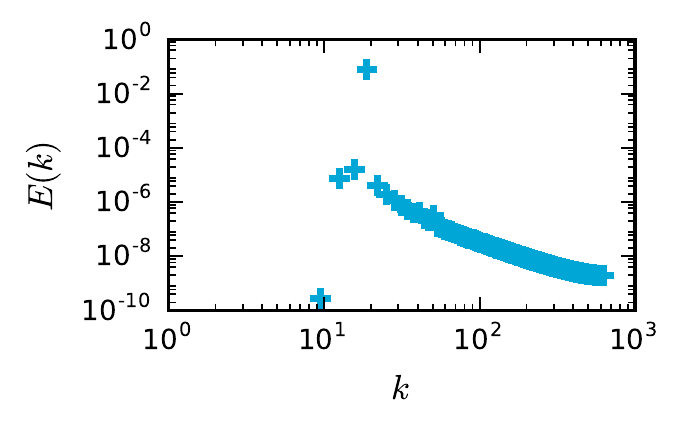}
		\caption{$t=0.0$}
	\end{subfigure}
	\begin{subfigure}{0.495\textwidth}
		\centering 
		\includegraphics{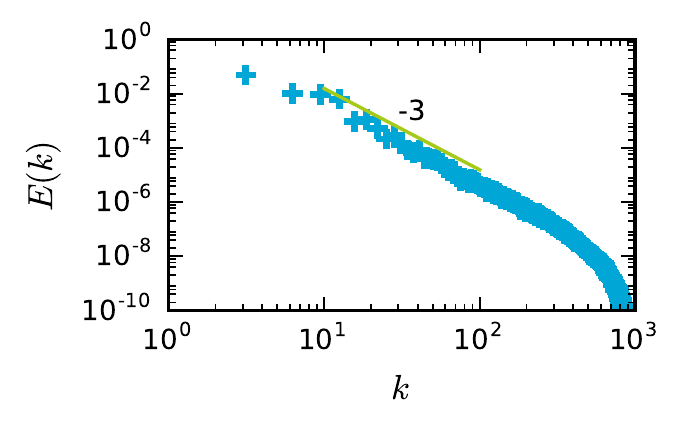}
		\caption{$t=10.0$} 
	\end{subfigure} \\
	\begin{subfigure}{0.495\textwidth}
		\centering 
		\includegraphics{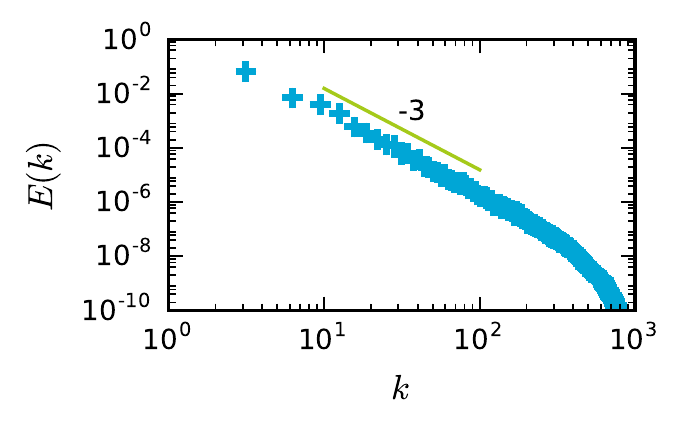}
		\caption{$t=15.0$} 
	\end{subfigure} 
	\begin{subfigure}{0.495\textwidth}
		\centering 
		\includegraphics{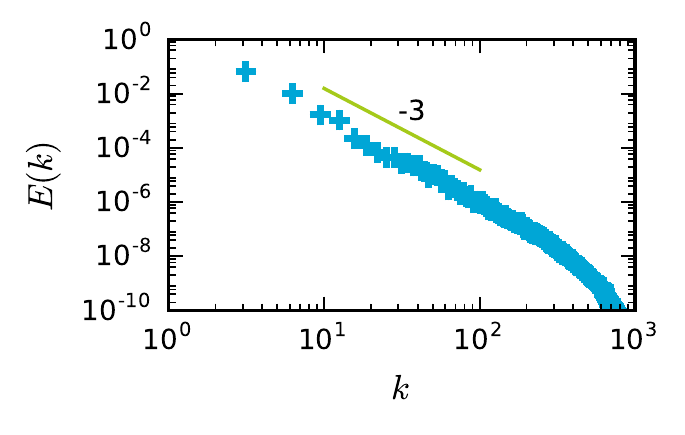}
		\caption{$t=20.0$} 
	\end{subfigure} \\
	\caption{Taylor-Green flow: One-dimensional energy spectra for Re=$2.5\cdot10^5$ at different time instants (blue) with a $k^{-3}$ slope (green).}
	\label{fig:taylor_green_lonu_energy}
\end{figure}

The energy density spectra for different time instants are plotted in Fig.~\ref{fig:taylor_green_lonu_energy}. At time instant 0, the energy density is concentrated at one wave-number. At later time instances, an energy cascade develops with a slope of approximately $k^{-3}$ covering almost two decades. 
This is in good correspondence with the $k^{-3}$ slope for the direct enstrophy cascade as predicted by Kraichnan \cite{Kraichnan1967}. A slight steepening of the spectrum is observed at later time instances. Corroboration for this behavior is found in the literature, where the coherent vortices are reported to destroy scale invariance and can produce steeper spectra than predicted by the $k^{-3}$ direct enstrophy cascade \cite{Clercx2009}. The increase in the energy density contained in the low wave numbers at the expense of the energy density contained in the higher wave number range, further illustrates the formation of large-scale coherent structures. 
\section{Conclusions and outlook} \label{sec: Conclusions Outlook}
This paper presented a particle-mesh scheme which allows for diffusion-free advection, satisfies mass and momentum conservation properties in a local (i.e. cellwise) sense, and allows the extension to high-order spatial accuracy. Central in obtaining this, is to formulate the particle-mesh projections in terms of a PDE-constrained minimization problem. The key idea in formulating the constraint is that from a mesh-perspective the particle motion must satisfy an advection operator. 
By expressing the control variable in terms of single-valued functions at cell interfaces, it was shown that the HDG method naturally provides the necessary ingredients for formulating the optimality system. Consistency and conservation of the constrained interaction were proven. An analysis of the resulting discrete optimality system, revealed that a specific choice for the Lagrange multiplier field renders the PDE-constrained optimization procedure independent of the time integration method, thus resulting in a particularly attractive and robust scheme. 

The PDE-constrained particle-mesh interaction is embedded in a Lagrangian-Eulerian operator splitting approach for the advection-diffusion and the incompressible Navier-Stokes equations. Consistency requirements specific to particle-mesh methods were formulated, and a particular fully-discrete formulation was formulated to obey these requirements.

A range of numerical experiments unveils the potential of the presented approach. Standard convergence rates in space and second-order convergence rates in time were confirmed for a range of examples for the linear advection-diffusion problem and the incompressible Navier-Stokes problem. The method turns out to be promising, albeit monotonicity is not guaranteed. A high Reynolds number test case gives further evidence of the robustness of the scheme.

The presented method opens many avenues for further investigation. Future work can include a more rigorous mathematical analysis of the method, finding alternative (high-order time-accurate) operator splitting approaches and/or time stepping strategies, and the reconstruction of monotonicity preserving background fields from scattered particle data by, e.g., exploiting the subgrid information available at the particle level. Future work will also investigate the applicability of the presented projection method for (mass conservative) interface tracking in, e.g., multiphase flows. 
\section*{Acknowledgements}
The Netherlands Organisation for Scientific Research (NWO) is gratefully acknowledged for their support through the JMBC-EM Graduate Programme research grant. The first author acknowledges Dr.~M.~M\"oller for the useful discussions on the subject, and the
4TU~Research~Centre~Fluid~\&~Solid~Mechanics 
for providing travel funding.

N.~Trask acknowledges support from the National Science Foundation MSPRF program, the Sandia~National~Laboratories LDRD program, and from the U.S. Department of Energy Office of Science, Office of Advanced Scientific Computing Research, Applied Mathematics program as part of the Collaboratory on Mathematics for Mesoscopic Modeling of Materials (CM4).
\appendix
\section{Consistency of the projection operators in the fully-discrete setting} \label{sec:appendix_consistency}
This appendix derives the consistency requirement for the projection operators $\mathcal{P}_E$ and $\mathcal{P}_L$ in the fully-discrete setting. To this end, we consider a repeated back-and-forth mapping of a scalar-valued quantity between the particles and the mesh, hence omitting the particle advection step (i.e. $\mathbf{a} = 0$) and the diffusion step (i.e. $\kappa = 0$). Moreover, we choose $l=0$. 
Consistency in the fully-discrete setting now requires the constraint equations to be inherently satisfied, i.e. $\xDiscreteScalar{\lambda} = 0$ everywhere. 

We assess this consistency criterion by evaluating the disrcete optimality system Eq.~\eqref{eq:discrete_optimality-adv-diff} at time level $n+1$. For a cell $K$ and with $\mathbf{a} = 0$ and $l=0$, these conditions simplify to
\begin{subequations} \label{eq:discrete_optimality_bare_bone}
\begin{align}\label{eq:costate-discrete_bare_bone}
\sum_{p \in \mathcal{S}_t^K} \xtDiscreteScalar{\psi}{n+1} (\mathbf{x}_p)  \xDiscreteScalar{w}(\mathbf{x}_p) 
+ 
\areaIntegral{K}{\frac{\xDiscreteScalar{w}}{\Delta t} \xtDiscreteScalar{\lambda}{n+1}  }
&=  
\sum_{p \in \mathcal{S}_t^K} \psi_p^{n}  \xDiscreteScalar{w}(\mathbf{x}_p),    \\
\areaIntegral{K}{\frac{\xtDiscreteScalar{\psi}{n+1}}{\Delta t}  \xDiscreteScalar{\tau}  }    &=  \areaIntegral{K}{  \frac{\xtDiscreteScalar{\psi}{*,n}}{\Delta t} \xDiscreteScalar{\tau} }, \label{eq:state-discrete_bare_bone}
\end{align}
\end{subequations}
where we used that $\xtDiscreteScalar{\bar{\psi}}{n+1}$ and $\xtDiscreteScalar{\psi}{n+1}$ coincide in a weak sense over the facets by virtue of the discrete optimality condition for $\mathbf{a}=\mathbf{0}$, Eq.~\eqref{eq:optimality-discrete-adv-diff}. \\
With the element contributions following from Eq.~\eqref{eq:discrete_optimality_bare_bone}, the system for a single element $K$ is casted as
	\begin{equation}
	\begin{bmatrix}
	\boldsymbol{M}_p 			& \boldsymbol{G} \\
	\boldsymbol{G}^\top  		& \boldsymbol{0}   			
	\end{bmatrix}
	\begin{bmatrix}
	\boldsymbol{\psi}^{n+1} \\ 
	\boldsymbol{\lambda}^{n+1} \\ 
	\end{bmatrix}
	=
	\begin{bmatrix}
	\boldsymbol{\chi}_p \boldsymbol{\psi}_p^{n} \\ 
	\boldsymbol{G}^\top \boldsymbol{\psi}^{*,n} \\
	\end{bmatrix}.
	\end{equation}
Performing a Gaussian elimination, results in the following algebraic form for $\boldsymbol{\lambda}^{n+1}$:
\begin{equation}
\boldsymbol{G}^\top \boldsymbol{M}_p^{-1}\boldsymbol{G} \boldsymbol{\lambda}^{n+1} = \boldsymbol{G}^\top \left(\boldsymbol{M}_p^{-1} \boldsymbol{\chi}_p \boldsymbol{\psi}_p^{n} - \boldsymbol{\psi}^{*,n} \right).
\end{equation}
The mass matrix, $\boldsymbol{M}_p$, is symmetric and non-negative definite \cite{Burgess1992,Brackbill1986}; thus, if it is non-singular, $\boldsymbol{G}^\top \boldsymbol{M}_p^{-1} \boldsymbol{G} \ne \boldsymbol{0}$. We obtain $\boldsymbol{\lambda}^{n+1} = \boldsymbol{0}$ if the term in the parentheses on the right-hand side equals $\boldsymbol{0}$. This, in turn, is achieved if and only if the mesh field $\boldsymbol{\psi}^{*,n}$ matches the (local $\ell^2$-) projection $\boldsymbol{M}_p^{-1} \boldsymbol{\chi}_p \boldsymbol{\psi}_p^{n}$. This condition indeed holds true by updating the particles via Eq.~\eqref{eq:discrete_particle_update} in conjunction with a definition of $\boldsymbol{\psi}^{*,n}$ via Eq.~\eqref{eq:consistency_term}.
\section*{References}
\bibliography{mybibfile}
\end{document}